\documentclass[10pt, a4paper,twoside]{article}
\usepackage{enumerate}
\usepackage{color}
\usepackage{color}
\usepackage{amsmath}
\usepackage[english]{babel}
\usepackage{amssymb,latexsym}
\usepackage{latexsym}
\usepackage{hyperref}

\def \bui#1#2{\mathrel{\mathop{\kern 0pt#1}\limits^{#2}}}

\newcommand{\R}{{\mathbb R}}
\renewcommand{\Re}{\mathrm{Re}}

\newtheorem{example}{Examples}[section]
\newtheorem{thm}{Theorem}[section]
\newtheorem{lemma}[thm]{Lemma}
\newtheorem{prop}[thm]{Proposition}

\newtheorem{cor}[thm]{Corollary}
\newtheorem{remark}[thm]{Remark}
\newtheorem{remarks}[thm]{Remarks}
\newtheorem{definition}[thm]{Definition}
\newtheorem{notation}[thm]{Notation}
\newtheorem{exabout:ample}[thm]{Example}

\title{Rigidity results for Riemannian spin$^c$ manifolds with foliated boundary}
\author{Fida El Chami\footnote{Lebanese University, Faculty of Sciences II, Department of Mathematics, P.O. Box 90656 Fanar-Matn, Lebanon,
E-mail: \texttt{fchami@ul.edu.lb}},\, Nicolas Ginoux\footnote{Institut \'Elie Cartan de Lorraine, Universit\'e de Lorraine, Site de Metz, B\^at. A, \^Ile du Saulcy, 57045 Metz Cedex 1, France, E-mail: \texttt{nicolas.ginoux@univ-lorraine.fr}},\, Georges Habib\footnote{Lebanese University, Faculty of Sciences II, Department of Mathematics, P.O. Box 90656 Fanar-Matn, Lebanon,
E-mail: \texttt{ghabib@ul.edu.lb}},\, Roger Nakad \footnote{Notre Dame University-Louaiz\'e, Faculty of Natural and Applied Sciences, Department of Mathematics and Statistics, P.O. Box 72, Zouk Mikael, Lebanon, E-mail: \texttt{rnakad@ndu.edu.lb}}}

\begin{document}
\date{}
\maketitle

\begin{abstract}
\noindent Given a Riemannian spin$^c$ manifold whose boundary is endowed with a Riemannian flow, we show that any solution of the basic Dirac equation  satisfies an integral inequality depending on geometric quantities, such as the mean curvature and the O'Neill tensor. We then characterize the equality case of the inequality when the ambient manifold is a domain of a K\"ahler-Einstein manifold or a Riemannian product of a K\"ahler-Einstein manifold with $\mathbb R$ (or with the circle $\mathbb S^1$).
\end{abstract}
{\bf Key words}: Manifolds with boundary, spin$^c$ structures, Riemannian flows, Basic Dirac equation, K\"ahler-Einstein manifolds, parallel spinors.

{\bf Mathematics Subject Classification:} 53C27, 53C12, 53C24.

\section{Introduction}

The spectral properties of the Dirac operator on a Riemannian spin manifold have led to several geometric and rigidity results. For example, spinorial techniques have been used to give simple proofs of some classical results, as the Alexandrov theorem \cite{HM3, HMneg}.   

On a compact spin manifold $N^{n+2}$ with boundary $M$ satisfying some curvature assumptions, O. Hijazi and S. Montiel \cite{HM1} proved that there exists a one-to-one correspondence between Killing spinors on $M$ and parallel spinors on $N.$ In particular, the boundary has to be connected and totally umbilical. This result has led to the following characterization of the round sphere as the boundary of the disk by: A complete Ricci-flat Riemannian manifold of dimension at least $3$ whose mean-convex boundary is isometric to the round sphere is a flat disc \cite[Cor. 6]{HM1}. In a more general setting, S. Raulot showed in \cite{Ra} that the correspondence occur between parallel spinors on $N$ and solutions of the Dirac equation on $M,$ i.e. a spinor field $\varphi$ satisfying $D_M \varphi=\frac{n+1}{2}H_0\varphi$ for some particular function $H_0$ with $H_0 \le H$ where  $H$ denotes the mean curvature of $M.$ As an application, he proved the following rigidity result \cite[Cor. 4]{Ra}: If $N$ has vanishing sectional curvature along the boundary (assumed to be simply connected), it has to be flat.

In a different geometric context, the authors in \cite{ElChamiGinoHabNaka16} established rigidity results for spin manifolds whose
boundary carries a Riemannian flow \cite{T}. That means the boundary is foliated by the integral curves of some unit vector field, say $\xi,$ in a way that the metric on $M$ stays constant along those curves. We will refer to the word {\it basic} to say objects that are constant along the curves (called also {\it leaves}). The idea is to check whether solutions of the so-called {\it basic Dirac equation} on the boundary, that is, a
spinor field $\varphi$ satisfying the equation
\begin{eqnarray}\label{D-E}
D_b\varphi=\frac{n+1}{2}H_0\varphi,
\end{eqnarray}
are in correspondence with parallel spinors on the whole manifold $N$. Here, $D_b$ is the basic Dirac operator (see \cite{GK1,GK2}) and $H_0$ is a basic function defined on the boundary. It turns out that, under some assumption relating the mean curvature to the O'Neill tensor of the flow \cite{ON}, this correspondence is valid. In particular, this characterizes the Riemannian product $\mathbb{S}^1\times \mathbb{S}^n$ as the boundary of $\mathbb{S}^1\times B,$ where $B$ is the closed unit ball in $\mathbb{R}^{n+1}$ \cite[Cor. 4.7]{ElChamiGinoHabNaka16}. 

\noindent The aim of this paper is to get similar rigidity results for manifolds endowed with spin$^c$ structures. Recall that those structures are the complex analogue to spin structures. They had real importance since the announcement of Seiberg-Witten theory \cite{F} (see references therein) whose applications to 4-dimensional geometry and topology  are  already  notorious.  From
an intrinsic point of view, spin, almost complex, complex, K\"ahler  and Sasaki manifolds have a
canonical spin$^c$ structure. From an extrinsic point of view, the restriction of spin$^c$ spinors is an
effective tool to study the geometry and the topology of submanifolds \cite{HMU06}. When shifting from spin  to spin$^c$ geometry, the situation becomes more general since the Dirac operator will not only depend on the geometry of the manifold but also on the connection of the auxiliary line bundle associated with the spin$^c$ structure. Thus, new examples appeared in several classification results \cite{{MoroianuCMP97}, GN}, less geometric and topological restrictions are imposed on the manifold and more flexibility is offered on the choice of a connection (and hence on the curvature $2$-form) on the auxiliary line bundle defining the spin$^c$ structure.

\noindent Throughout this paper, we will consider spin$^c$ manifolds with foliated boundary. We will look at the solutions of the basic Dirac equation on the boundary, i.e.,  Equation (\ref{D-E}), where $D_b$ is now the spin$^c$ basic Dirac operator. First,
we start by restricting the spin$^c$ structure on $N$ to the boundary $M$ and define a spin$^c$ structure on the normal bundle $Q$ of the flow by taking the same auxiliary bundle as the one on $M$. Then, we define a connection $1$-form on $Q$ by modifying the one on $M$ in a way to get a basic connection.
It turns out that the parameter chosen to make the connection one-form basic is related to the curvature $F^M$ of the auxiliary bundle of $M$ via Equation \eqref{eq:AQbasique}.
Second, we show that solutions of the  basic Dirac equation satisfy a {\it spin$^c$ integral inequality} derived from  \cite[Prop. 9]{HM}. Indeed, if we denote by $F^N$ the curvature $2$-form on the auxiliary line bundle defining the spin$^c$ structure on  $N,$ we have:
 
\begin{thm}\label{pro:estimate-int} Let $(N^{n+2},g)$ be a compact spin$^c$ manifold with boundary.
Assume that the scalar curvature ${\rm Scal}^N $ of $N$ satisfies ${\rm Scal}^N\geq c_{n+2} |F^N|$ ($c_{n+2} = 2[\frac{n+2}{2}]^{\frac 12}$), that the boundary $M$ of $N$ carries a minimal Riemannian flow given by a unit vector field $\xi$ and that the mean curvature $H$ is positive.
Assume also that the connection on $Q$ is basic, i.e.,  \eqref{eq:AQbasique} is fulfilled for some basic function $\theta$ as well as the existence of a basic section $\varphi$ of $\Sigma M$ such that $D_b\varphi=\frac{(n+1)H_0}{2}\varphi$ for some non-negative basic function $H_0$ on $M$.
Then the following inequality holds:
\begin{equation}\label{eq:integraline22-int}
0\leq\int_M\frac{1}{H}\cdot\Big((n+1)^2\left(H_0^2-H^2\right)|\varphi|^2+|(\Omega\cdot_M-i\theta)\varphi|^2
\Big)dv_g,\end{equation}
where $\Omega$ is the $2$-form associated to the O'Neill tensor field. 
\end{thm}

We mention that spin manifolds correspond to the case where the function $\theta$ vanishes everywhere and therefore Inequality \eqref{eq:integraline22-int} is the same as in \cite[Thm. 3.1]{ElChamiGinoHabNaka16}. However, the equality case of Inequality \eqref{eq:integraline22-int} is characterized by the existence of two parallel spinors on $N$ that project (through some orthogonal pointwise projection, see Proposition \ref{pro:hijazimontiel}) to the solution $\varphi.$ Therefore and according to A.~Moroianu classification of spin$^c$ manifolds with parallel spinors \cite{MoroianuCMP97}, it turns out that, besides spin manifolds, two natural categories of manifolds occur in the limiting case: Either $N$ is a domain in some K\"ahler-Einstein manifold of nonnegative scalar curvature ($n$ is even) or in some Riemannian product of a K\"ahler-Einstein manifold  with $\mathbb R$ ($n$ is odd) (see Proposition \ref{pro:caracterisation}). 

\noindent In the even dimensional case, i.e. $N$ is a K\"ahler manifold endowed with its canonical spin$^c$ structure, the limiting case of (\ref{eq:integraline22-int}) is characterized in the following theorem:

\begin{thm}\label{p:eqcaseKaehler}
Under the assumptions of {\rm Theorem \ref{pro:estimate-int}}, let furthermore $(N,g,J)$ be K\"ahler, endowed with its canonical spin$^c$ structure.
Then $(N,g)$ is Einstein with nonnegative scalar curvature and Inequality {\rm(\ref{eq:integraline22-int})} is satisfied by $\varphi$.
Moreover, {\rm(\ref{eq:integraline22-int})} is an equality for some nonzero $\varphi$ if and only if the following properties hold:
\begin{enumerate}
\item Up to changing $\xi$ into $-\xi$ and hence $\theta$ into $-\theta$, we have $\xi=-J\nu$.
In particular, $[J,A_{|_Q}]=0$ and the vector field $\xi$ is a pointwise eigenvector for the Weingarten map $A$, that is, $A\xi=\lambda\xi$ for some $\lambda\in C^\infty(M,\mathbb{R})$.
\item The function $\theta$ is constant and equal to $\lambda+\frac{1}{2}\mathrm{tr}\left(A_{|_Q}\right)$.
\end{enumerate}
If those properties are fulfilled, then the spinor field $\varphi$ is the restriction of a parallel spinor of $N$ and is itself transversally parallel, in particular $H_0=0$.
Furthermore, $\xi(\lambda)=\xi\left(\mathrm{tr}\left(A_{|_Q}\right)\right)=0$.
\end{thm}

\noindent In particular when $N$ is a domain in either $\mathbb{C}\mathrm{P}^{\frac{n+2}{2}}$ or $\mathbb{C}\mathrm{H}^{\frac{n+2}{2}},$ then the equality holds in (\ref{eq:integraline22-int}) if and only if $M$ is a tube around a totally geodesic $\mathbb{C}\mathrm{P}^k$ ($0\leq k\leq\frac{n-2}{2}$) or a tube around a totally geodesic $\mathbb{C}\mathrm{H}^k$ (where $0\leq k\leq\frac{n}{2}),$ see Proposition \ref{c:classifequalcaseKaehlercstholseccurv}. On the other hand and under some assumption relating  $H$, $H_0$, $\Omega$ and $\theta$, we show that the equality in \eqref{eq:integraline22-int} is realized and in this case, the flow is transversally Einstein-K\"ahler and the manifold $M$ is $\eta$-Einstein Sasakian manifold (see Corollary \ref{coro-even}).

 \noindent In the odd dimensional case, i.e. when $N$ is a domain in the Riemannian product $N_1 \times \mathbb R$ or $N_1 \times \mathbb S^1$ where $N_1$  is a closed K\"ahler-Einstein manifold, we obtain the following result:

\begin{thm} \label{pro:caslimiteimpair}
Under the assumptions of Theorem \ref{pro:estimate-int}, assume furthermore $N$ to be a domain in the Riemannian product $N_1\times \R$ or $N_1\times \mathbb{S}^1$ where $N_1$ is a closed K\"ahler-Einstein manifold.
Then the equality case is realized in \eqref{eq:integraline22-int} if and only if $N$ is isometric to the Riemannian product $\Delta\times \mathbb{S}^1$ with $\xi=\pm\partial t$ where $\Delta$ is a domain in $N_1$.
If those conditions are satisfied, then $\varphi$ is some pointwise projection of a parallel spinor on $N_1$, $\theta=0$ and $H_0=H$.
\end{thm}

\noindent Note in particular that equality can only occur if $N$ lies in $N_1\times\mathbb{S}^1$, the case where $N\subset N_1\times\R$ leading to the noncompactness of $M$, see proof of Lemma \ref{lem:cara1}.
The proof of Theorem \ref{pro:caslimiteimpair} is splitted into several steps. First, we prove that $N$ is isometric to the product $\Delta\times \mathbb{S}^1$ for some compact domain $\Delta$ with boundary $M_1$ (see Lemma \ref{lem:cara1}). Second, we show after some technical computations that the vector field $\xi$ defines a Riemannian flow on $M_1$ by some vector field $\xi_1.$ In this case, the normal bundle $Q_1$ (of even rank) of the flow is a subbundle of $Q.$ After restricting the spin$^c$ structure on $N$ to the bundle $Q_1$, we prove that the spinor $\varphi$ defines another spinor field $\varphi_1$ which is a solution of the basic Dirac equation on $M_1$ and realizes the equality case in (\ref{eq:integraline22-int}) for the even dimensional case (see Lemma \ref{lem:caract12}). This last part is proved with the help of Theorem \ref{p:eqcaseKaehler}.
Finally, we show that this leads to $\xi=\pm \partial t$ and $H=H_0.$ 

\section{Preliminaries on Riemannian flows and manifolds with boundary}\label{sec:prel}
\noindent In this section, we recall some preliminaries on spin$^c$ Riemannian flows (see \cite{T}, \cite{GH}) and the geometry of manifolds with boundary. For more details, we refer to \cite{LM}, \cite{F}, \cite{G} and \cite{BHMMM}.
\subsection{spin$^c$ Riemannian flow}\label{sub2.1}
\noindent Let $(M^{n+1},g)$ be a Riemannian manifold endowed with a Riemannian flow given by a unit vector field $\xi$.
That is, the Lie derivative of the metric with respect to $\xi$ satisfies  $\mathcal{L}_{\xi} g_{|_{\xi^\perp}}=0$ (see \cite{Rei}).
It is now well-known that, if that condition is fulfilled, then there exists a metric connection $\nabla$ on the normal bundle $Q=\xi^\perp$ which is called transversal Levi-Civita connection and which is defined, for any section $Y \in \Gamma(Q)$, by

\begin{equation}\label{eq:connexiontransversal}
\nabla _{X} Y :=
\left\{\begin{array}{ll}
\pi [X,Y]&\textrm {if $X=\xi$},\\\\
\pi (\nabla_{X}^{M}Y)&\textrm{if  $X\perp \xi$},
\end{array}\right.
\end{equation}
where $\pi:TM\rightarrow Q$ denotes the orthogonal projection \cite{T}.
The transversal Levi-Civita connection is related to the usual Levi-Civita connection via the following Gauss-type formulas \cite{GH}: for all sections $Z,W$ in $\Gamma(Q),$ we have
\begin{equation}\label{eq:Oneillformula}
\left\{\begin{array}{ll}
\nabla^M_Z W=\nabla_Z W-g(h(Z),W)\xi, &\textrm {}\\\\
\nabla^M_\xi Z=\nabla_\xi Z+h(Z)-\kappa(Z)\xi,&\textrm {}
\end{array}\right.
\end{equation}
where $h:=\nabla^M\xi$ is the O'Neill tensor and $\kappa:=h(\xi)$ is the mean curvature of the flow.

\noindent From now on, we assume that $M$ is a spin$^c$ manifold.
That means, there exist a ${\rm Spin}_{n+1}^c$-principal bundle ${\rm Spin}^c(M)$ and a $\mathbb{U}(1)$-principle bundle $P_{\mathbb{U}_1}M$ over $M$ (called the {\it auxiliary line bundle} of the spin$^c$ structure) together with a double covering $\eta:{\rm Spin}^c(M)\rightarrow {\rm SO}M\times P_{\mathbb{U}_1}M$ such that $\eta(ua) = \eta(u)\eta_0(a),$ for every $u\in {\rm Spin}^c M$ and $a\in {\rm Spin}_{n+1}^c$ where $\eta_0: {\rm Spin}_{n+1}^c\rightarrow {\rm SO}_{n+1}\times \mathbb{S}^1$ is the 2-fold covering.
Here, ${\rm SO}M$ denotes the ${\rm SO}_{n+1}$-principal bundle of orthonormal direct frames on $TM.$
As for the spin case, the decomposition of the tangent bundle of $M$ into $TM=\R\xi\oplus Q$ allows to induce the ${\rm Spin}^c$  structure on $M$ to a ${\rm Spin}^c$ structure on the normal bundle $Q$ (see \cite{B,ElChamiGinoHabNaka16} for details).
The spin$^c$ structure on $Q$ is given by the pull-back of the one on $M$ via the inclusion map ${\rm SO}Q\rightarrow {\rm SO}M,$ where ${\rm SO}Q$ denotes the ${\rm SO}_n$-principal bundle of orthonormal direct frames on $Q.$ The auxiliary line bundle of the spin$^c$ structure of $Q$ is chosen to be the same one as on $M$.\\

{\bf Choice of connections:} We choose a connection $1$-form $A^M$ on the auxiliary line bundle $P_{\mathbb{U}_1}M\to M$ and define a connection $1$-form on the auxiliary line bundle of the spin$^c$ structure of $Q$ by \emph{modifying} the connection $1$-form on the auxiliary line bundle of the spin$^c$ structure of $M$: Pick any real $1$-form on $M$, say $\alpha$, and define a new connection $1$-form on $P_{\mathbb{U}_1}M\to M$ by
\begin{equation}\label{eq:connexionmodifee}
A^Q:=A^M-i\alpha.
\end{equation}
This makes sense since the difference of any two connections $1$-forms on a $\mathbb{U}_1$-bundle is given by an imaginary-valued form on the base $M$.
Since, we need to have an $\mathcal{F}$-bundle $P_{\mathbb{U}_1}M$ in the sense of \cite{ElKacimiGmira97,Elkacimi90} (in order to define later the basic Dirac operator), we wish the connection $A^Q$ to be \emph{basic}, which is equivalent to the condition $\xi\lrcorner F^Q=0$ on $M$, where $F^Q\in\Omega^2(M,i\mathbb{R})$ is the curvature $2$-form associated to $A^Q$, see e.g. \cite[p.328]{ElKacimiGmira97}.
This last identity means that
\[\xi\lrcorner F^M=i\xi\lrcorner d\alpha.\]
For the sake of simplicity, we \emph{assume} $\alpha$ to be proportional to $\xi^\flat$, that is, that $\alpha=\alpha(\xi)\cdot\xi^\flat$.
Setting $\theta:=\alpha(\xi)\in C^\infty(M,\mathbb{R})$, the condition to be fulfilled for $A^Q$ to be basic becomes equivalent to
\begin{equation}\label{eq:AQbasique}
\xi\lrcorner F^M=-id\theta+i\theta\kappa^\flat+i\xi(\theta)\xi^\flat
\end{equation}
on $M$.
Namely, $\xi\lrcorner F^M(\xi)=0=\left(-id\theta+i\theta\kappa^\flat+i\xi(\theta)\xi^\flat\right)(\xi)$ (we have $g(\kappa,\xi)=0$ because of $g(\xi,\xi)=1$) and, for any $X\in \Gamma(Q)$,
\begin{eqnarray*}
\left(i\xi\lrcorner d\alpha\right)(X)&=&id\alpha(\xi,X)\\
&=&i\left\{\xi(\underbrace{\alpha(X)}_{0})-X\left(\alpha(\xi)\right)-\alpha\left([\xi,X]\right)\right\}\\
&=&i\left\{-X(\theta)-g([\xi,X],\xi)\alpha(\xi)\right\}\\
&=&i\left\{-X(\theta)+g(X,\nabla_\xi^M\xi)\theta\right\},
\end{eqnarray*}
which implies (\ref{eq:AQbasique}).\\

Coming back to spin$^c$ structures, the spinor bundle $\Sigma M$ is canonically identified with the spinor bundle of $Q,$ denoted by $\Sigma Q,$
for $n$ even and with the direct sum $\Sigma Q\oplus \Sigma Q$ for
$n$ odd. In the same way, one can also identify the Clifford multiplications ``$\cdot_M$'' in $\Sigma M$
and ``$\cdot_Q$'' in $\Sigma Q$ as follows: For any section
$Z\in \Gamma(Q)$ and $\varphi\in \Gamma(\Sigma M),$ we have
\begin{equation*}
\left\{\begin{array}{ll}
Z\cdot_M \varphi=Z \cdot_Q\varphi, &\textrm {for $n$ even}\\\\
\xi \cdot_M Z\cdot_M \varphi=(Z\cdot_Q\oplus-Z\cdot_Q)\varphi,&\textrm {for $n$ odd}.
\end{array}\right.
\end{equation*}

\noindent The Levi-Civita connections on $\Sigma M$ and $\Sigma Q$ satisfy the formulas \cite[eq. (2.4.7)]{GH}
\begin{equation}\label{eq:oneillspin}
\left\{\begin{array}{ll}\nabla_\xi^M\varphi&=
\nabla_\xi\varphi+\frac{1}{2}\Omega\cdot_M\varphi+\frac{1}{2}\xi\cdot_M\kappa\cdot_M\varphi
+\frac{i \theta}{2} \varphi\\
 &\\
\nabla_Z^M\varphi&= \nabla_Z\varphi+\frac{1}{2}\xi\cdot_M
h(Z)\cdot_M\varphi,\end{array}\right.
\end{equation}
where $\Omega$ is the $2$-form associated to the tensor $h$ defined for all $Y,Z\in \Gamma(Q)$ by $\Omega(Y,Z)=g(h(Y),Z).$

\noindent Next, we define the basic Dirac operator (see \cite{GK1} and \cite{GK2}) as being
$$D_b=\sum_{i=1}^n e_i\cdot_Q\nabla_{e_i}-\frac{1}{2}\kappa\cdot_Q,$$
where $\{e_i\}_{i=1,\cdots,n}$ is a local orthonormal frame of
$\Gamma(Q).$
Here and in the following, we assume the mean curvature to be basic (otherwise, we might work with its basic projection \cite{RicPark}), i.e. $\nabla_\xi\kappa=0.$ Recall that the basic Dirac operator $D_b$ is defined on the set of basic sections (sections of the spinor bundle $\Sigma Q$ satisfying $\nabla_\xi\varphi=0$) and it preserves that set.
It is also a transversally elliptic and essentially self-adjoint operator, if $M$ is compact.
Therefore, it has a discrete spectrum by the spectral theory of transversal elliptic operators \cite{ElKacimiGmira97,Elkacimi90}.

\noindent As a direct consequence of Equations \eqref{eq:oneillspin}, the transverse Levi-Civita connection commutes with the Clifford action of $\xi.$
In particular, this allows to prove the following identities for the basic Dirac operator.
For $n$ even (resp. $n$ odd) and for any basic spinor field $\varphi$, we have
\begin{equation}\label{eq:commute}
D_b(\xi\cdot_M \varphi)=-\xi\cdot_M D_b \varphi\,\, \,\,({\rm resp.}\,\, D_b(\xi\cdot_M \varphi)=\xi\cdot_M D_b \varphi).
\end{equation}
Finally, we recall the relation existing between the Dirac
operator on $M$ and the basic Dirac operator:
\begin{equation}\label{eq:relationdirac}
\left\{\begin{array}{ll}
D_M=D_b-\frac{1}{2}\xi\cdot_M\Omega\cdot_M +\frac{i\theta}{2} \xi\cdot_M, &\textrm {for $n$ even}\\\\
D_M=-\xi\cdot_M(D_b\oplus -D_b)-\frac{1}{2}\xi\cdot_M\Omega\cdot_M
+\frac{i\theta}{2} \xi\cdot_M, &\textrm
{for $n$ odd}.
\end{array}\right.
\end{equation}

\subsection{Manifolds with boundary}\label{manifbound}

\noindent We review some well-known facts about ${\rm Spin}^c$ manifolds with boundary (see \cite{HM1,HM2,HM3} for the spin case).
Let $(N^{n+2},g)$ be a Riemannian ${\rm Spin}^c$
manifold of dimension $n+2$ with smooth boundary $M=\partial N.$
As before, the existence of the (inward) unit vector field $\nu$ normal to the boundary allows to define a ${\rm Spin}^c$ structure on $M$ by taking the pull back.
We can define two spinor bundles on the boundary, the intrinsic bundle $\Sigma M$ and the extrinsic one  $\textbf{S}=\Sigma N_{|_{M}}.$  The data of the extrinsic bundle is related to the one on $N$ by:
\begin{eqnarray}\label{eq:gauss}
X\cdot_\textbf{S} \varphi &=&X\cdot \nu\cdot\varphi\nonumber\\
\nabla^N_X\varphi &=&\nabla^\textbf{S}_X \varphi +\frac{1}{2}A(X)\cdot_\textbf{S}\varphi\\
\textbf{D}_{\textbf{S}}\varphi&=&\frac{n+1}{2}H\varphi-\nu\cdot D_N\varphi-\nabla^N_\nu\varphi \nonumber,
\end{eqnarray}
where $``\cdot"$ is the Clifford multiplication on $N$, the tensor $A$ is the Weingarten map given for all $X\in \Gamma(TM)$ by $A(X)=-\nabla^N_X\nu,$ the spinor field $\varphi$ is a section in $\textbf{S}$ and $H=\frac{1}{n+1}{\rm Trace}(A)$ is the mean curvature of $M.$ The operator $\textbf{D}_\textbf{S},$ called the extrinsic Dirac operator, acts on sections on $\textbf{S}$ as $\textbf{D}_\textbf{S}=\sum_{i=1}^{n+1} e_i\cdot_\textbf{S}\nabla^{\textbf{S}}_{e_i},$ where $\{e_1,\cdots,e_{n+1}\}$ is a local orthonomal frame of $TM.$

\noindent On the other hand, the extrinsic spinor bundle can be identified with the intrinsic one in a canonical way depending on the dimension of $N.$ Namely, if $n$ is odd, the tuple $(\textbf{S},``\cdot_\textbf{S}",\nabla^\textbf{S},\textbf{D}_\textbf{S})$ can be identified with $(\Sigma M,``\cdot_{M}",\nabla^{M},D_{M})$ whereas for $n$ even it can be identified with $(\Sigma M\oplus \Sigma M,``\cdot_{M}\oplus -\cdot_{M}",\nabla^{M}\oplus \nabla^{M},D_{M}\oplus -D_{M}).$ Moreover, using the first two equations in \eqref{eq:gauss} and the Gauss formula, one can prove that the following relations hold for all $X,Y\in \Gamma(TM)$
\begin{equation} \label{eq:compa}
\left\{\begin{array}{ll}
\nabla^\textbf{S}_X(Y\cdot)=\nabla^M_X Y\cdot+Y\cdot\nabla^\textbf{S}_X,\\\\
\nabla^\textbf{S}_X(\nu\cdot)=\nu\cdot\nabla_X^\textbf{S}
\end{array}\right.
\end{equation}
and that,
\begin{equation}\label{eq:espacepropre}
\textbf{D}_\textbf{S}(\nu\cdot)=-\nu\cdot\textbf{D}_\textbf{S}.
\end{equation}
Equality \eqref{eq:espacepropre} means that the spectrum of $\textbf{D}_\textbf{S}$ is symmetric with respect to zero and if $n$ is even the Dirac operator on $M$ commutes with the action of $\nu$, that is,
\begin{eqnarray} \label{eq:comu}
D_M(\nu\cdot\Phi)=\nu\cdot D_M\Phi
\end{eqnarray}
for any spinor field $\Phi\in \Gamma(\Sigma M).$

\noindent We define the operators $P_{\pm}$ as being the pointwise orthogonal projections from $\textbf{S}$ onto the $\pm 1$-eigenspaces corresponding to the $\pm 1$-eigenvalues of the operator $i\nu\cdot$ on $\textbf{S},$ i.e. $ P_{\pm}:=\frac{1}{2}\left(\mathrm{Id}\pm i\nu\cdot\right).$
They satisfy $$P_{\pm}(X \cdot)=X\cdot P_{\mp} \ \ \ \text{and} \ \ \ \ P_\pm(\nu\cdot)=\nu\cdot P_\pm,$$
for all $X\in \Gamma(TM)$.
This implies that $\textbf{D}_\textbf{S}P_{\pm}=P_{\mp}\textbf{D}_\textbf{S}.$

\section{Integral inequality on manifolds with boundary}\label{sec:integralinequalitymfdswithboundary}

\noindent In \cite{HM},  O.~Hijazi and S.~Montiel prove an integral inequality relating the Dirac operator on the boundary $M$ of a spin manifold $N$ applied to a spinor field to the norm of that spinor. In the following, we will state a similar inequality for spin$^c$ manifolds.

\begin{prop}\label{pro:hijazimontiel} Let $(N^{n+2},g)$ be a compact spin$^c$ manifold such that ${\rm Scal}^N\geq c_{n+2}|F^N|$, where $F^N$ is the curvature of the auxiliary line bundle and $c_{n+2}:=2[\frac{n+2}{2}]^\frac{1}{2}.$
If the mean curvature $H$ of the boundary $M$ is positive, then any spinor field $\varphi\in\Gamma({\rm {\bf {S}}})$ satisfies the inequality
\begin{equation}\label{hij-mon}
0\leq \int_M\Big(\frac{1}{H}|{\rm\bf D}_{\rm\bf{S}}\varphi|^2-\frac{(n+1)^2}{4}H|\varphi|^2\Big)dv_g,
\end{equation}
where $dv_g$ is the volume element on $M.$
Moreover, the equality holds in {\rm(\ref{hij-mon})} if and only if there exist two parallel spinor fields $\psi,\vartheta\in
\Gamma(\Sigma N)$ such that $P_+\varphi=P_+\psi$ and $P_-\varphi=P_-\vartheta$ on the boundary.
In that case, the scalar curvature of $N$ is equal to $c_{n+2}|F^N|$, in particular is nonnegative.
\end{prop}

{\it Proof:} The proof of the inequality is based on the spinorial Reilly formula established in \cite[p.142]{Nakad11} and an appropriate boundary value problem.
We will prove the inequality for $n$ odd and will use the different identifications according to Subsection \ref{manifbound}.
The same can be done for $n$ even.\\
For any spinor field $\psi\in\Gamma(\Sigma N)$, one has \cite[p.142]{Nakad11}
\begin{eqnarray}\label{eq:Reillyformula}
\int_M(\langle D_M\psi,\psi\rangle-\frac{n+1}{2}H |\psi|^2)s_g \geq &\ &\frac{1}{4}\int_N ({\rm Scal}^N-c_{n+2}|F^N|)|\psi|^2 v_g\nonumber \\ &&-\frac{n+1}{n+2}\int_N |D_N\psi|^2 v_g.
\end{eqnarray}
Equality holds in \eqref{eq:Reillyformula} if and only if $\psi$ is a twistor spinor and
\begin{equation}\label{eq:curvatureterm}
F^N\cdot\psi=-\frac{c_{n+2}}{2}|F^N|\psi.
\end{equation}
Recall here that a twistor spinor $\psi$ is a section of the spinor bundle of $N$ satisfying the differential equation $\nabla^N_X\psi=-\frac{1}{n+2}X\cdot D_N\psi$ for any $X\in\Gamma(TN).$ In the following, we will follow the same proof as in \cite[p.11]{HM}. For this, consider for any spinor field $\varphi\in\Gamma({\rm {\bf {S}}})$, the following boundary value problem:
$$
\left\{\begin{array}{ll}
D_N\psi=0 &\textrm {on $N$},\\\\
P_+\psi=P_+\varphi &\textrm{on  $M$}.
\end{array}\right.
$$
The uniqueness and the smoothness of the solution of the boundary problem is shown e.g. in \cite[Prop. 6]{HM}.
By inserting the solution $\psi$ into Inequality \eqref{eq:Reillyformula}, we get after using ${\rm Scal}^N\geq c_{n+2}|F^N|$ that
\begin{equation}\label{eq:integralinequality}
\int_M(\langle D_M\psi,\psi\rangle-\frac{n+1}{2}H |\psi|^2)s_g \geq 0.
\end{equation}
Here, we notice that the equality in \eqref{eq:integralinequality} is realized if and only if the spinor $\psi$ is parallel.
In that case, Equality \eqref{eq:curvatureterm} is automatically satisfied as a consequence of the Schr\"odinger-Lichnerowicz formula.
Now by decomposing the spinor field $\psi$ into $\psi=P_+\psi+P_-\psi$ and using the pointwise inequality $0\leq |\frac{1}{\sqrt{\frac{n+1}{2}H}}D_MP_+\psi-\sqrt{\frac{n+1}{2}H}P_-\psi|^2,$ Inequality \eqref{eq:integralinequality} becomes
\begin{equation} \label{eq:inegaliteplus}
\int_M (\frac{1}{H}|D_M P_+\varphi|^2-\frac{(n+1)^2}{4}H |P_+\varphi|^2)s_g \geq 0.
\end{equation}
Finally, considering a boundary problem for $D_N$ where we replace the condition boundary $P_+$ by $P_-,$ we get a similar equation as \eqref{eq:inegaliteplus} with the minus sign.
Summing the two equations, we deduce the desired inequality.
\hfill$\square$ \\

\noindent As we have seen, the limiting case of \eqref{eq:integralinequality} is characterized by the existence of parallel spinors on the ambient manifold.
Recall that the boundaryless complete Riemannian spin$^c$ manifolds carrying parallel spinors were classified by A.~Moroianu in \cite[Thm. 1.1]{MoroianuCMP97}: The universal cover of such a manifold is isometric to the Riemannian product of a simply connected K\"ahler manifold with a simply connected spin manifold carrying parallel spinors; the classification remains true even locally.
In the following, we will determine which family of such products satisfies the condition ${\rm Scal}^N\geq c_{n+2}|F^N|$ required for Proposition \ref{pro:hijazimontiel}.
We begin with examining that condition on K\"ahler manifolds:

\begin{prop}\label{prop:einsteinkahler} Let $(N^{n+2},g,J)$ be a K\"ahler manifold endowed with its canonical spin$^c$ structure, in particular $F^N=-i{\rm Ric}_N\circ J$ is the complexified Ricci form of $N.$
If ${\rm Scal}^N\geq c_{n+2}|F^N|$ holds, where $c_{n+2}:=2^[\frac{n+2}{2}]^\frac{1}{2}$, then that inequality is an equality and $N$ is Einstein with nonnegative scalar curvature.
\end{prop}

{\it Proof:} \noindent Choose a local o.n.b. $\left(e_j\right)_{1\leq j\leq n+2}$ of $TN$ made out of pointwise eigenvectors for the Ricci tensor of $N$, that is, $\mathrm{Ric}_N(e_j)=\mu_j\cdot e_j$ for some real eigenvalue $\mu_j$ and all $j\in\{1,\ldots,n+2\}$.
Then $\mathrm{Scal}_N=\sum_{j=1}^{n+2}\mu_j$ and we have
\begin{eqnarray*}
c_{n+2}\cdot|F^N|&=&2\cdot\sqrt{\frac{n+2}{2}}\cdot\sqrt{\frac{1}{2}\sum_{j=1}^{n+2}\left|\mathrm{Ric}_N(Je_j)\right|^2}\\
&=&\sqrt{n+2}\cdot\sqrt{\sum_{j=1}^{n+2}\mu_j^2},
\end{eqnarray*}
so that the assumption $\mathrm{Scal}_N\geq c_{n+2}|F^N|$ becomes
\[\sum_{j=1}^{n+2}\mu_j\geq\sqrt{n+2}\cdot\sqrt{\sum_{j=1}^{n+2}\mu_j^2}.\]
But, by the Cauchy-Schwarz inequality, $\sum_{j=1}^{n+2}\mu_j\leq\sqrt{n+2}\cdot\sqrt{\sum_{j=1}^{n+2}\mu_j^2}$.
Therefore, we necessarily have the equality $\sum_{j=1}^{n+2}\mu_j=\sqrt{n+2}\cdot\sqrt{\sum_{j=1}^{n+2}\mu_j^2}$ and therefore all $\mu_j$ are equal and nonnegative.
This means precisely that the Einstein condition is satisfied and that the Einstein constant is nonnegative.
\hfill$\square$\\

\noindent Using this result, we deduce:

\begin{prop} \label{pro:caracterisation} Let $(N^{m},g)$ be a Riemannian spin$^c$ manifold carrying a parallel spinor such that ${\rm Scal}^N\geq c_{m}|F^N|$ where $F^N$ is the curvature of the auxiliary line bundle of the spin$^c$ structure.
Then $N$ is locally isometric to either a spin manifold with parallel spinors, a K\"ahler-Einstein manifold of nonnegative scalar curvature or the Riemannian product of a K\"ahler-Einstein manifold of nonnegative scalar curvature with $\mathbb{R}.$
\end{prop}

{\it Proof:} Since all requirements are of algebraic type, we may assume that $N$ is simply connected.
As mentioned before, we know that $N$ is locally isometric to the Riemannian product $N_1\times N_2$ where $N_1$ is K\"ahler and $N_2$ is a spin manifold carrying parallel spinors (and which is in particular Ricci flat).
Assume that $N_1$ is not a point, otherwise we are reduced to the spin case.
Then the condition $\mathrm{Scal}_N\geq c_{m}|F^N|$ can be written in terms of the data of $N_1$.
Namely, $\mathrm{Scal}_N=\mathrm{Scal}_{N_1}$ and $|F^N|=|F^{N_1}|.$
Hence if we denote by $n_1$ the dimension of $N_1$, we get
$$\mathrm{Scal}_{N_1}\geq c_{m}|F^{N_1}|\geq c_{n_1}|F^{N_1}|.$$
\noindent Since $N_1$ is K\"ahler, we deduce from Proposition \ref{prop:einsteinkahler} that $N_1$ is Einstein-K\"ahler and also the fact that $c_m=c_{n_1}.$ Mainly, that means either $m=n_1$ or $m=n_1+1$.
Therefore, $N$ is locally K\"ahler or the product of a K\"ahler manifold with $\mathbb{R}$.
\hfill$\square$

\section{Hypersurfaces of K\"ahler-Einstein manifolds} \label{sec:kahlereinstein}

\noindent In this section, we consider real hypersurfaces in any K\"ahler manifold.
We first characterize the condition for the naturally induced flow-structure to be Riemannian and minimal and then study the $\eta$-umbilicity condition on the hypersurface.


\begin{prop}\label{p:HopfEinstein}
Let $M^{n+1}$ be any immersed real oriented hypersurface in a K\"ahler manifold $(N^{n+2},g,J)$.
Denote by $\nu$ the unit normal inducing the orientation along $M$, by $\xi$ the tangent vector field $-J\nu$ on $M$ and by $Q:=\xi^\perp\subset TM$ the horizontal distribution on $M$.
Then the following holds:
\begin{enumerate}
\item For every $X\in TM$, we have $\nabla_X^M\xi=JAX-g(A\xi,X)\nu$.
In particular, $(M,g,\xi)$ is a minimal Riemannian flow if and only if $\left[J,A_{|_Q}\right]=0$ and $M$ is Hopf, that is, that there is a function $\lambda$ on $M$ such that $A\xi=\lambda\xi$ on $M$, where $A:=-\nabla^N\nu$ is the Weingarten map of $M$.
\item If $(M,g,\xi)$ is a minimal Riemannian flow, then the complex structure $J$ induces a K\"ahler structure on the normal bundle $Q$ of the flow with respect to the transversal Levi-Civita connection defined by \eqref{eq:connexiontransversal}.
\item If $(M,g,\xi)$ is a minimal Riemannian flow and $N$ is Einstein, then $\xi\left(\mathrm{tr}\left(A_{|_Q}\right)\right)=0$.
\end{enumerate}
\end{prop}

{\it Proof}:
\begin{enumerate}
\item First note that $\xi$ has unit length, in particular $h:=\nabla^M\xi$ induces an endomorphism field of $Q$.
For every $X\in \Gamma(TM)$,
\[\nabla_X^M\xi=\nabla_X^N\xi-g(AX,\xi)\nu=-J(\nabla_X^N\nu)-g(A\xi,X)\nu=JAX-g(A\xi,X)\nu.\]
In particular, $\nabla_\xi^M\xi=0$ if and only if $A\xi$ is pointwise proportional to $\xi$, that is, if $M$ is Hopf.
Furthermore, $h$ is pointwise skew-symmetric if and only if $(JA)^*=-JA$ on $Q$, where $(\cdot)^*$ denotes the pointwise $g$-adjoint; but $(JA)^*=A^*J^*=-AJ$ on $Q$, which implies the first statement.
\item Since $\xi=-J\nu,$ the complex structure $J$ maps the normal bundle $Q$ to itself.
To prove that $J$ defines a K\"ahler structure, it is sufficient to show that it is parallel with respect to the transversal connection $\nabla.$
But for any section $X\in\Gamma(Q)$,
\begin{eqnarray*}
(\nabla_\xi J)(X)&=&\nabla_\xi (J(X))-J(\nabla_\xi X)\\
&\bui{=}{\eqref{eq:Oneillformula}}& \nabla^M_\xi (J(X))-h(J(X))-J(\nabla_\xi^M X)+J(h(X))\\
&=& \nabla^N_\xi (J(X))-g(A(\xi),J(X))\nu-J(\nabla_\xi^N X)-g(A(\xi),X)\xi\\
&=& (\nabla^N_\xi J)(X)-\lambda g(\xi,J(X))\nu-\lambda g(\xi,X)\xi\\
&=&0.
\end{eqnarray*}
The next to last equality comes from the Gauss formula and the first statement.
Now, for $Y,Z\in \Gamma(Q),$ we compute in the same way
\begin{eqnarray*}
(\nabla_Y J)(Z)&=&\nabla_Y (J(Z))-J(\nabla_Y Z)\\
&\bui{=}{\eqref{eq:Oneillformula}}& \nabla^M_Y J(Z)+g(h(Y),J(Z))\xi-J(\nabla_Y^M Z)-g(h(Y),Z)\nu\\
&=& \nabla^N_Y J(Z)-g(A(Y),J(Z))\nu+g(h(Y),J(Z))\xi-J(\nabla_Y^N Z)\\&&-g(A(Y),Z)\xi-g(h(Y),Z)\nu\\
&=&0,
\end{eqnarray*}
from which $\nabla J=0$ follows.
\item We make use of the following well-known formula, valid as soon as the ambient manifold is Einstein:
\begin{equation}\label{eq:divA}
\delta A=-(n+1)dH,
\end{equation}
where $\delta A:=-\sum_{j=1}^{n+1}(\nabla_{e_j}^MA)(e_j)$ for any local o.n.b. $(e_j)_{1\leq j\leq n+1}$ of $TM$.
Choosing this local o.n.b. such that $e_{n+1}=\xi$, we have
\begin{eqnarray*}
\delta(A\xi)&=&-\sum_{j=1}^{n+1}g(\nabla_{e_j}^MA\xi,e_j)\\
&=&-\sum_{j=1}^{n+1}g((\nabla_{e_j}^MA)(\xi),e_j)-\sum_{j=1}^{n+1}g(A(\nabla_{e_j}^M\xi),e_j)\\
&=&-\sum_{j=1}^{n+1}g((\nabla_{e_j}^MA)(e_j),\xi)-g(A(\nabla_\xi^M\xi),\xi)-\sum_{j=1}^{n}g(A(he_j),e_j)\\
&=&(\delta A)(\xi)-\lambda \underbrace{g(\nabla_\xi^M\xi,\xi)}_{0}-\underbrace{g\left(A_{|_Q},h\right)}_{0}\\
&=&(\delta A)(\xi)\\
&=&-(n+1)\xi(H)\qquad\textrm{ by }(\ref{eq:divA})\\
&=&-\xi(\lambda)-\xi\left(\mathrm{tr}\left(A_{|_Q}\right)\right).
\end{eqnarray*}
On the other hand,
\begin{eqnarray*}
\delta(A\xi)&=&\delta(\lambda\xi)\\
&=&-\xi(\lambda)+\lambda\cdot\delta\xi\\
&=&-\xi(\lambda)-\lambda\cdot\sum_{j=1}^{n+1}g(\nabla_{e_j}^M\xi,e_j)\\
&=&-\xi(\lambda)-\lambda\cdot\sum_{j=1}^n\underbrace{g(he_j,e_j)}_{0}-\lambda\cdot\underbrace{g(\nabla_{\xi}^M\xi,\xi)}_{0}\\
&=&-\xi(\lambda).
\end{eqnarray*}
Comparing both identities, we obtain
\[\xi(\lambda)=\xi(\lambda)+\xi\left(\mathrm{tr}\left(A_{|_Q}\right)\right),\]
from which $\xi\left(\mathrm{tr}\left(A_{|_Q}\right)\right)=0$ follows.
\end{enumerate}
\hfill$\square$

Next we look at conditions for certain real hypersurfaces in K\"ahler manifolds to have constant principal curvatures:

\begin{prop}\label{p:hypersisopar}
Under the assumptions of {\rm Proposition \ref{p:HopfEinstein}}, let $N$ be Einstein, the flow induced by $\xi:=-J\nu$ be a minimal Riemannian flow and $\theta$ be constant and equal to $\lambda+\frac{1}{2}\mathrm{tr}\left(A_{|_Q}\right)$.
If $A_{|_Q}=\mu\cdot\mathrm{Id}_Q$ for some $\mu\in C^\infty(M,\mathbb{R})$, then unless $n=2$ both $\lambda$ and $\mu$ are constant.
In the case where $N$ has constant holomorphic sectional curvature, the same conclusion holds even in the case where $n=2$.
\end{prop}

{\it Proof}: We compute $\delta A$ using (\ref{eq:divA}) and choosing $(e_j)_{1\leq j\leq n+1}$ with $Ae_j=\mu\cdot e_j$ for all $1\leq j\leq n$ and $e_{n+1}=\xi$, so that $Ae_{n+1}=\lambda\cdot e_{n+1}$.
Note that, for any $1\leq j\leq n$, the (local) vector field $\nabla_{e_j}^Me_j$ is transverse: for $g(\nabla_{e_j}^Me_j,\xi)=-g(\nabla_{e_j}^M\xi,e_j)=0$ since $h$ is skew-symmetric.
We compute, also using the fact that $\kappa=0$:
\begin{eqnarray*}
\delta A&=&-\sum_{j=1}^n\left(\nabla_{e_j}^M(\mu\cdot e_j)-A(\nabla_{e_j}^Me_j)\right)-\nabla_\xi^M(\lambda\xi)+A(\nabla_\xi^M\xi)\\
&=&-\sum_{j=1}^n\left(e_j(\mu)e_j+\mu\nabla_{e_j}^Me_j-\mu\nabla_{e_j}^Me_j\right)-\xi(\lambda)\xi\\
&=&-d\mu_{|_Q}-\xi(\lambda)\xi.
\end{eqnarray*}
Comparing with (\ref{eq:divA}), we obtain the following identity:
\[nd\mu+d\lambda=d\mu_{|_Q}+\xi(\lambda)\xi.\]
Putting $X=\xi$, we first deduce that $\xi(\mu)=0$, so that the last identity becomes
\[(n-1)d\mu+d\lambda-\xi(\lambda)\xi=0.\]
But, since we have assumed $\theta=\lambda+\frac{n\mu}{2}$ to be constant, we also have $\frac{n}{2}d\mu+d\lambda=0$, from which we deduce that $\xi(\lambda)=0$ and $(n-2)d\mu=0$.
If $n\neq2$, we can conclude that $d\mu=0=d\lambda$.\\
In the case where $N$ has constant holomorphic sectional curvature $4c$ with $c\geq0$, its curvature tensor is given by the following identity \cite[Theorem 1.1]{NiebergallRyan97}:
\begin{equation}\label{eq:curvtensorcstholsect}R_{X,Y}^NZ=-4c\left(X\wedge Y + JX\wedge JY + 2g(JX,Y)J\right)Z,\end{equation}
where $(X\wedge Y)Z=g(X,Z)Y-g(Y,Z)X$.
Now we can compute $(\nabla_X^MA)(Y)-(\nabla_Y^MA)(X)$ and use the Codazzi identity: for any $X,Y\in\Gamma(Q)$,
\begin{eqnarray*}
(\nabla_X^MA)(Y)&=&\nabla_X^M(AY)-A(\nabla_X^MY)\\
&=&\nabla_X^M(\mu Y)-A(\nabla_X^MY)\qquad\textrm{ where }\nabla_X^MY=\nabla_XY-g(hX,Y)\xi\\
&=&X(\mu)Y+\mu(\nabla_XY-g(hX,Y)\xi)-A(\nabla_XY-g(hX,Y)\xi)\qquad\\
&=&X(\mu)Y+\mu(\nabla_XY-\mu g(JX,Y)\xi)-\mu\nabla_XY+\lambda\mu g(JX,Y)\xi\\
&=&X(\mu)Y+\mu(\lambda-\mu)g(JX,Y)\xi,
\end{eqnarray*}
where we used that $h=JA$ and $A\xi=\lambda\xi.$ Hence, we deduce that for any $X,Y\in\Gamma(Q)$,
\[(\nabla_X^MA)(Y)-(\nabla_Y^MA)(X)=X(\mu)Y-Y(\mu)X+2\mu(\lambda-\mu)g(JX,Y)\xi.\]
Now Codazzi formula states that $(\nabla_X^MA)(Y)-(\nabla_Y^MA)(X)=-R_{X,Y}^N\nu$, so that (\ref{eq:curvtensorcstholsect}) yields, still for any $X,Y\in\Gamma(Q)$,
$$X(\mu)Y-Y(\mu)X+2\mu(\lambda-\mu)g(JX,Y)\xi=-R^N_{X,Y}\nu=-8cg(JX,Y)\xi,$$
from which $X(\mu)=0$ follows (choosing $X,Y\in \Gamma(Q)$ pointwise linearly independent, which is always possible if $n\geq2$).
We can conclude that $\mu$ and hence also $\lambda$ are constant on $M$.
\hfill$\square$

The condition $\theta$ be constant and equal to $\lambda+\frac{1}{2}\mathrm{tr}\left(A_{|_Q}\right)$ comes from Proposition \ref{p:eqcaseKaehlersuite} below.

\begin{prop}\label{p:bordSasaki}
Under the assumptions of {\rm Proposition \ref{p:HopfEinstein}}, let the flow induced by $\xi:=-J\nu$ be a minimal Riemannian flow.
If $A_{|_Q}=\mu\cdot\mathrm{Id}_Q$ for some $\mu\in\mathbb{R}^\times$, then $\left(M,\frac{1}{\mu^2}\cdot g_\xi\oplus g_Q,\mu\xi\right)$ is a Sasakian manifold.
\end{prop}

{\it Proof}: Recall that a Riemannian flow is Sasaki if and only if it is minimal and $h=\nabla^M\xi$ is a transverse K\"ahler structure on $Q$, that is, $h^2=-\mathrm{Id}_Q$ and $\nabla h=0$.
Recall also that, for any $t\in\mathbb{R}^\times$, the triple $(M,g_t:=t^2g_\xi\oplus g_Q,\frac{1}{t}\xi)$ is a Riemannian flow with $h_t=th$, $\kappa_t=\kappa$ and $\nabla^t=\nabla$ on $Q$.
In our situation, if $A_{|_Q}=\mu\cdot\mathrm{Id}_Q$ for some $\mu\in\mathbb{R}^\times$, then $h=JA_{|_Q}=\mu J$ with $\nabla J=0$, so that $\frac{1}{\mu}h=J$ defines a transverse K\"ahler structure on $Q$.
This concludes the proof of Proposition \ref{p:bordSasaki}.
\hfill$\square$

\section{spin$^c$ manifolds with foliated boundary}\label{sec:spincmfdswithfoliatedboundary}

In this section, we consider a compact Riemannian spin$^c$ manifold $(N^{n+2},g)$ with nonempty boundary $M=\partial N$ and assume that $M$ is endowed with a Riemannian flow induced by a unit vector field $\xi$ on $M.$
After restricting the spin$^c$ structure to the normal bundle of the flow, we will consider solutions of the basic Dirac equation and will show that it satisfies an integral inequality coming from \eqref{eq:integralinequality}.
Then, we will study the limiting case of this inequality and characterize the geometry of the manifold $N$ and its boundary $M$ based on the results in Section \ref{sec:kahlereinstein}.

\noindent We have seen that the limiting case of Inequality \eqref{eq:integralinequality} is characterized by the existence of parallel spinors on the whole manifold $N.$
In view of Proposition \ref{pro:caracterisation}, we will consider two families: K\"ahler-Einstein domains and domains in products of K\"ahler-Einstein manifolds with $\R$ or $\mathbb{S}^1$.
Therefore, we will split our study into two cases: The even- and odd-dimensional case.
The even-dimensional case will correspond to the first family of manifolds while the odd-dimensional one will correspond to the second.

\subsection{The even-dimensional case}\label{subsec:evendimcase}

As mentioned before, we will consider in this subsection the case where the ambient manifold is K\"ahler-Einstein and is endowed with its canonical spin$^c$ structure.
First, we need to characterize the condition \eqref{eq:AQbasique} for an $\mathcal{F}$-bundle in this setting.

\begin{prop}\label{p:eqcaseKaehlersuite}
Let $M^{n+1}$ be any connected immersed hypersurface in a K\"ahler-Einstein manifold $(N^{n+2},g,J).$
Let $N$ carry its canonical spin$^c$ structure, $M$ carry the spin$^c$ structure induced by the inner unit normal $\nu$ and $Q:=(J\nu)^\perp\cap TM$ in turn carry the spin$^c$ structure induced by $\xi:=-J\nu$.
Assume $(M,g,\xi)$ to be a minimal Riemannian flow.
Choose a connection $1$-form on the line bundle $K_N^{-1}{}_{|_M}=P_{\mathbb{U}_1}M\to M$ associated to the spin$^c$ structure of $Q$ with the help of a basic function $\theta\in C^\infty(M,\mathbb{R})$ as given in  \eqref{eq:connexionmodifee}.

Then the following holds:
\begin{enumerate}
\item The induced connection $1$-form on $P_{\mathbb{U}_1}M$ is basic if and only if $\theta$ is constant on $M$.
\item The connection $1$-form $A^Q$ given in  \eqref{eq:connexionmodifee} is induced by the transverse Levi-Civita connection $\nabla$ of $Q$ if and only if $\theta=\lambda+\frac{1}{2}\mathrm{tr}\left(A_{|_Q}\right)$, where $A\xi=\lambda\xi$ on $M$.
This is also equivalent to the restriction of a parallel section $\psi\in\Gamma(\Sigma_0N)$ to $M$ being basic, which is also equivalent to it being transversally parallel.
\item If $\theta=\lambda+\frac{1}{2}\mathrm{tr}\left(A_{|_Q}\right)$ and is constant on $M$, then the transversal Ricci tensor of the flow satisfies ${\rm Ric}^\nabla={\rm Ric}^N+2\theta A.$  In particular, if $\theta=0$ or $A_{|_Q}=\mu\mathrm{Id}_Q$ for some $\mu\in\mathbb{R}$, then the flow is transversally Einstein-K\"ahler.
\item If $\theta=\lambda+\frac{1}{2}\mathrm{tr}\left(A_{|_Q}\right)$ and is constant on $M,$ the curvature of the line bundle of $Q$ is equal to $-i{\rm Ric}^\nabla\circ J.$ Mainly, that means the induced spin$^c$ structure on $Q$ from $N$ is the canonical spin$^c$ structure on $Q$ (the one induced by $J$).
\end{enumerate}
\end{prop}

{\it Proof}: The statements follow from elementary computations.
Note before computing anything that the existence of a Riemannian flow on $M$ implies that complex structure $J$ is transversally K\"ahler as shown in Proposition \ref{p:HopfEinstein}.
\begin{enumerate}
\item Since the line bundle of the canonical spin$^c$ structure on $N$ is $K_N^{-1}$ and the connection is the Levi-Civita one, its curvature form is given by $F^N=-i\mathrm{Ric}_N\circ J$ and hence $\xi\lrcorner F^M=\xi\lrcorner F^N=-i\mathrm{Ric}^N(\nu)=0$ along $M$ because of the Einstein condition on $N$.
By (\ref{eq:AQbasique}), the connection $1$-form on $Q$ is basic if and only if $0=-d\theta+\xi(\theta)\xi^\flat$, that is, $d\theta=0$ since $\theta$ is assumed to be basic here.
\item First note that we can express $A^N$ and then $A^Q$ with the help of a local Hermitian o.n.b. of $TN$.
Namely let $\langle\cdot\,,\cdot\rangle:=g(\cdot\,,\cdot)-ig(J\cdot\,,\cdot)$ be the Hermitian metric associated to $g$, let $(e_1,\ldots,e_{\frac{n+2}{2}})$ be any local $\langle\cdot\,,\cdot\rangle$-o.n.b. of $TN$, or equivalently, such that $(e_1,Je_1,\ldots,e_{\frac{n+2}{2}},Je_{\frac{n+2}{2}})$ is a local $g$-o.n.b. of $TN$, then $s:=e_1^*\wedge\ldots\wedge e_{\frac{n+2}{2}}^*$ is a local section of $K_N^{-1}$ and hence, for any $X\in \Gamma(TN)$,
\begin{eqnarray*}
\nabla_X^Ns&=&\left(\sum_{j=1}^{\frac{n+2}{2}}\langle\nabla_X^Ne_j^*,e_j^*\rangle\right)e_1^*\wedge\ldots\wedge e_{\frac{n+2}{2}}^*\\
&=&-\left(\sum_{j=1}^{\frac{n+2}{2}}\langle\nabla_X^Ne_j,e_j\rangle\right)e_1^*\wedge\ldots\wedge e_{\frac{n+2}{2}}^*\\
&=&i\cdot\left(\sum_{j=1}^{\frac{n+2}{2}}g(\nabla_X^N(Je_j),e_j)\right)e_1^*\wedge\ldots\wedge e_{\frac{n+2}{2}}^*,
\end{eqnarray*}
so that $-iA^N(s_*X)=\sum_{j=1}^{\frac{n+2}{2}}g(\nabla_X^N(Je_j),e_j)$.
Restricting to $M$ and choosing for instance $e_{\frac{n+2}{2}}=\xi$ and hence $Je_{\frac{n+2}{2}}=\nu$, we obtain, still picking a local section $s$ of $K_N^{-1}{}_{|_M}$,
\[-iA^M(s_*X)=\sum_{j=1}^{\frac{n}{2}}g(\nabla_X^N(Je_j),e_j)+g(\nabla_X^N\nu,\xi)=\sum_{j=1}^{\frac{n}{2}}g(\nabla_X^M(Je_j),e_j)-g(A\xi,X).\]
As for $A^Q$, we first have, for any $X\in \Gamma(Q)$,
\[-iA^Q(s_*X)=-iA^M(s_*X)=\sum_{j=1}^{\frac{n}{2}}g(\nabla_X^M(Je_j),e_j)-\lambda\underbrace{g(\xi,X)}_{0}=\sum_{j=1}^{\frac{n}{2}}g(\nabla_X(Je_j),e_j).\]
For $X=\xi$, we have
\begin{eqnarray*}
-iA^Q(s_*\xi)&=&-iA^M(s_*\xi)-\alpha(\xi)\\
&=&\sum_{j=1}^{\frac{n}{2}}g(\nabla_\xi^M(Je_j),e_j)-g(A\xi,\xi)-\theta\\
&=&\sum_{j=1}^{\frac{n}{2}}g(\nabla_\xi(Je_j),e_j)+g(h(Je_j),e_j)-\lambda-\theta\qquad\textrm{with }h=JA=AJ\\
&=&\sum_{j=1}^{\frac{n}{2}}g(\nabla_\xi(Je_j),e_j)-g(Ae_j,e_j)-\lambda-\theta\\
&=&\sum_{j=1}^{\frac{n}{2}}g(\nabla_\xi(Je_j),e_j)-\frac{1}{2}\mathrm{tr}\left(A_{|_Q}\right)-\lambda-\theta.
\end{eqnarray*}
This shows that $A^Q$ is associated to $\nabla$ if and only if $\theta=\lambda+\frac{1}{2}\mathrm{tr}\left(A_{|_Q}\right)$.\\
As for the second equivalence, note that $\psi\in\Gamma(\Sigma_0Q)$: because of $\psi\in\Gamma(\Sigma_0N)$ and hence $\xi\cdot\nu\cdot\psi=-i\psi$, we have
\[\Omega^Q\cdot\psi=\Omega^N\cdot\psi-\xi\cdot\nu\cdot\psi=-i\frac{(n+2)}{2}\psi+i\psi=-i\frac{n}{2}\psi,\]
where $\Omega^Q$ and $\Omega^N$ are the $2$-forms associated to $J$ on $Q$ and $N$ respectively.\\
For any $X\in \Gamma(Q)$, using that $JY\cdot\psi=iY\cdot\psi$ for any $Y\in \Gamma(TN)$ because of $\psi\in\Gamma(\Sigma_0N)$, we write
\begin{eqnarray*}
\nabla_X\psi&=&\nabla_X^M\psi-\frac{1}{2}\xi\cdot h(X)\cdot\psi\qquad\textrm{with }h=JA\\
&=&\underbrace{\nabla_X^N\psi}_{0}-\frac{1}{2}A(X)\cdot\nu\cdot\psi-\frac{1}{2}\xi\cdot JA(X)\cdot\psi\qquad\textrm{with }\nu=J\xi\\
&=&-\frac{i}{2}A(X)\cdot\xi\cdot\psi-\frac{i}{2}\xi\cdot AX\cdot\psi\\
&=&0\qquad\textrm{since }A\xi=\lambda\xi,
\end{eqnarray*}
so that $\nabla_X\psi=0$.
In particular, $\psi$ is transversally parallel if and only if $\nabla_\xi\psi=0$, that is, if and only if $\psi$ is basic.
On the other hand, we can make use of the following formula to relate $\nabla_\xi\psi$ with $\nabla_\xi^N\psi$:
\begin{eqnarray*}\nabla_\xi\psi&=&\underbrace{\nabla_\xi^N\psi}_{0}-\frac{1}{2}A\xi\cdot\nu\cdot\psi-\frac{1}{2}\Omega\cdot\psi-\frac{i\theta}{2}\psi\qquad\textrm{ with }A\xi=\lambda\xi\\
&=&\frac{i(\lambda-\theta)}{2}\psi-\frac{1}{2}\Omega\cdot\psi,
\end{eqnarray*}
so that $\nabla_\xi\psi=0$ if and only if $\Omega\cdot\psi=i(\lambda-\theta)\psi$.
But we can compute the action of $\Omega$ in another way, using $h=JA_{|_Q}$ as well as $\psi\in\Gamma(\Sigma_0N)$: choose $\left(e_j\right)_{1\leq j\leq n}$ to be any local orthonormal basis of $Q$ made out of eigenvectors for $A$, where $Ae_j=\mu_j e_j$, then
\begin{eqnarray*}
\Omega\cdot_M\psi&=&\frac{1}{2}\sum_{j=1}^ne_j\cdot_M he_j\cdot_M\psi\\
&=&\frac{1}{2}\sum_{j=1}^ne_j\cdot_M JAe_j\cdot_M\psi\\
&=&\frac{1}{2}\sum_{j=1}^n\mu_j e_j\cdot_M Je_j\cdot_M\psi\qquad\textrm{with }Je_j\cdot_M\psi=ie_j\cdot_M\psi\\
&=&-\frac{i}{2}\sum_{j=1}^n\mu_j\psi\\
&=&-\frac{i}{2}\left(\mathrm{tr}(A)-g(A\xi,\xi)\right)\psi\\
&=&\frac{i}{2}\left(\lambda-(n+1)H\right)\psi,
\end{eqnarray*}
so that $\frac{i}{2}\left(\lambda-(n+1)H\right)=i(\lambda-\theta)$, which yields $\lambda-2\theta=-(n+1)H$ and this is equivalent to $\theta=\lambda+\frac{1}{2}\mathrm{tr}\left(A_{|_Q}\right)$.
\item Let $\psi$ be any parallel spinor on $N$ which is then a transversal parallel spinor on $Q.$ Using the spin$^c$ Ricci identity on the normal bundle, we write for any $Y \in \Gamma(Q)$
\begin{eqnarray*}
\mathrm{Ric}^\nabla Y \cdot_Q \psi &=& (Y \lrcorner F^Q) \cdot_Q \psi =  \sum_{k=1}^n F^Q(Y,e_k) e_k \cdot_Q \psi\\
&=&  \sum_{k=1}^nF^N (Y,e_k) e_k \cdot_Q \psi-i\theta d\xi(Y,e_k)e_k \cdot_Q \psi\\
&=&\sum_{k=1}^nF^N (Y,e_k) e_k \cdot_Q \psi-2i\theta h(Y)\cdot_Q \psi.
\end{eqnarray*}
Therefore by the fact that $(Y\lrcorner F^N)(\xi)=(Y\lrcorner F^N)(\nu)=0$ and since $\psi$ is parallel on $N$, we deduce that
\begin{eqnarray*}
\mathrm{Ric}^\nabla Y \cdot \psi &=&{\rm Ric}^N Y \cdot \psi-2\theta Jh(Y)\cdot \psi\\
&=&{\rm Ric}^N Y \cdot \psi+2\theta A(Y)\cdot \psi,
\end{eqnarray*}
and the proof follows.
\item Using the previous computation for the Ricci curvature, one can evaluate the curvature of the auxiliary line bundle of $Q$ as follows
\begin{eqnarray*}
F^Q(Y,Z)&=& F^N(Y,Z)+2i\theta g(h(Z),Y)\\
&=&-i{\rm Ric}^N(JY,Z)+2i\theta g(h(Z),Y)\\
&=&-i{\rm Ric}^\nabla(JY,Z)+2i\theta g(A(JY),Z)+2i\theta g(h(Z),Y)\\
&=&-i{\rm Ric}^\nabla(JY,Z).
\end{eqnarray*}
This ends the proof of the proposition.
\hfill$\square$
\end{enumerate}

\noindent Before we state our main estimate, we need the following lemma:
\begin{lemma} \label{lem:commu} Let $(N^{n+2},g)$ be a spin$^c$ manifold. Assume that $n$ is even and the boundary $M$ of $N$ carries a minimal Riemannian flow given by a unit vector field $\xi$.Then, for any basic spinor field $\varphi,$ one has
$$D_M(\xi\cdot\varphi)=-\xi\cdot D_b\varphi+\frac{1}{2}\nu\cdot\Omega\cdot\varphi
-\dfrac{i\theta}{2} \nu \cdot \varphi.$$
\end{lemma}

{\it Proof}: Since the manifold $N$ is spin$^c$, the manifold $M$ and the normal bundle $Q$ of
the flow are also spin$^c$. Because $n$ is even, the spinor
bundle of $Q$ is identified with the one of $M,$ which is also
identified with one subbundle of $\Sigma N,$ say $\Sigma^+N$.
Therefore we can think of any $\varphi\in \Gamma(\Sigma Q)$ as a
section in one subbundle of $\textbf{S}$, say $\textbf{S}^+.$ We compute
\begin{eqnarray*}
D_M(\xi\cdot\varphi)&=&D_M(\nu\cdot\xi\cdot\nu\cdot\varphi)
=D_M(\nu\cdot(\xi\cdot_M\varphi))\nonumber\\
&\bui{=}{\eqref{eq:comu}}& \nu\cdot D_M(\xi\cdot_M\varphi)\nonumber\\
&\bui{=}{\eqref{eq:relationdirac}}&\nu\cdot (  D_b(\xi\cdot_M\varphi)-\frac{1}{2}\xi\cdot_M\Omega\cdot_M(\xi\cdot_M\varphi) +\dfrac{i\theta}{2} \xi \cdot_M (\xi \cdot_M \varphi))\nonumber\\
&\bui{=}{\eqref{eq:commute}}&\nu\cdot (-\xi\cdot_M D_b\varphi+\frac{1}{2}\Omega\cdot_M\varphi-\dfrac{i\theta}{2} \varphi)\nonumber\\
&=&-\xi \cdot D_b\varphi+\frac{1}{2}\nu\cdot\Omega\cdot\varphi
-\dfrac{i\theta}{2} \nu \cdot \varphi.
\end{eqnarray*}
This finishes the proof of the lemma.
\hfill$\square$\\

\noindent Now, we can prove Theorem \ref{pro:estimate-int}  for $n$ even.


{\it Proof of Theorem \ref{pro:estimate-int} for $n$ even}: As in Lemma \ref{lem:commu}, we will think of any $\varphi\in \Gamma(\Sigma Q)$ as a
section in the subbundle $\textbf{S}^+$ of $\textbf{S}$.
Using Equation \eqref{eq:relationdirac}, we compute
\begin{equation*}
 {\textbf D}_{\textbf{S}} \varphi=D_M\varphi=\frac{n+1}{2}H_0\ \varphi-\frac{1}{2}\xi\cdot_M\Omega\cdot_M\varphi
+ \frac{i \theta}{2} \xi \cdot_M \varphi.
\end{equation*}
Hence, by taking the norm of  ${\textbf D}_{\textbf{S}}\varphi,$ we get
\begin{eqnarray*}
\vert {\textbf D}_{\textbf{S}} \varphi\vert^2 &=&
\frac{(n+1)^2}{4} H_0^2 \vert \varphi\vert^2 +\frac 14 \vert
\Omega\cdot_M\varphi\vert^2 + \frac{\theta^2}{4} \vert
\varphi\vert^2-\dfrac{n+1}{2} H_0 \langle \varphi, \xi\cdot_M \Omega \cdot_M
\varphi \rangle\nonumber\\&& + \dfrac{n+1}{2} H_0 \theta  \langle \varphi, i
\xi\cdot_M \varphi \rangle - \frac{\theta}{2} \langle \Omega
\cdot_M\varphi, i\varphi \rangle.
\end{eqnarray*}
On the other hand, using Lemma \ref{lem:commu} we have
\begin{eqnarray*}
{\textbf D}_{\textbf{S}} (\xi\cdot\varphi)=-D_M (\xi\cdot\varphi) &=& \xi \cdot D_b \varphi -\frac 12 \nu
\cdot \Omega\cdot\varphi + \frac i2 \theta \nu \cdot\varphi \\
&=& \frac{n+1}{2} H_0 \xi \cdot\varphi - \frac 12 \nu \cdot
\Omega\cdot\varphi + \frac {i}{2} \theta \nu \cdot\varphi.
\end{eqnarray*}
As before, the norm is equal to
\begin{eqnarray*}
\vert {\textbf D}_{\textbf{S}} (\xi\cdot\varphi)\vert^2 &=&
\frac{(n+1)^2}{4} H_0^2 \vert \varphi\vert^2 +\frac 14 \vert
\Omega\cdot_M\varphi\vert^2 + \frac{\theta^2}{4} \vert
\varphi\vert^2+\dfrac{n+1}{2} H_0 \langle \varphi, \xi\cdot_M \Omega \cdot_M
\varphi \rangle\nonumber\\&& - \dfrac{n+1}{2} H_0 \theta  \langle \varphi, i
\xi\cdot_M \varphi \rangle - \frac{\theta}{2} \langle \Omega
\cdot_M \varphi, i\varphi \rangle.
\end{eqnarray*}
Thus, by applying Inequality \eqref{eq:integralinequality} to the spinors $\varphi$ and $\xi\cdot\varphi,$ we find the desired result after summing both inequalities.
\hfill$\square$

{\it Proof of Theorem \ref{p:eqcaseKaehler}}:
We know from Proposition \ref{pro:hijazimontiel} that if (\ref{eq:integraline22-int}) is an equality, there exist parallel spinor fields $\psi,\phi,\Psi,\Phi$ on $N$ such that $P_+\varphi=P_+\psi$, $P_-\varphi=P_-\phi$, $P_+(\xi\cdot\varphi)=P_+\Psi$ and $P_-(\xi\cdot\varphi)=P_-\Phi$ on $M$.
Note that $\psi\in\Gamma(\Sigma_0N)$ and that all other parallel spinors are multiple scalars of $\psi$, since the manifold $N$ is assumed to be K\"ahler.
But $P_+\psi=\frac{1}{2}(\psi+i\nu\cdot\psi)$; since $\psi\in\Gamma(\Sigma_+N)$, we have $i\nu\cdot\psi\in\Gamma(\Sigma_-N)$.
The same holds for $\varphi$ (always true for $n$ even).
We can thus deduce from $P_+\varphi=P_+\psi$ that $\varphi=\psi$ on $M.$\\
Of course, we also have $\psi=\phi$ from $P_-\varphi=P_-\phi$.
The identity $P_+(\xi\cdot\varphi)=P_+\Psi$ yields in the same way $\xi\cdot\varphi=i\nu\cdot\Psi$ and hence $i\nu\cdot\xi\cdot\varphi=\Psi$.
Letting $\Psi=b\psi$ for some $b\in\mathbb{C}$, we have $b=\pm1$ because of $(i\nu\cdot\xi\cdot)^2=1$.
Up to changing $\xi$ into $-\xi$ (and hence $\theta$ into $-\theta$), we can assume that $i\nu\cdot\xi\cdot\varphi=-\psi=-\varphi$.
Notice now that, if $X$ and $Y$ are two real vectors in $T_xN$ for some point $x\in N$ with $(X+iY)\cdot\psi=0$ for some nonzero $\psi\in\Sigma_0N_x$, then $Y=JX$: just combine the identity with $(X+iJX)\cdot\psi=0$, which holds true because of $\psi\in\Sigma_0N$.
Therefore, $i\xi\cdot\nu\cdot\varphi=\varphi$ being equivalent to $(\nu-i\xi)\cdot\psi=0$, we can conclude that $\xi=-J\nu$.
The last identity $P_-(\xi\cdot\varphi)=P_-\Phi$ does not bring any new information.
Since by assumption $(M,g,\xi=-J\nu)$ is a minimal Riemannian flow, Proposition \ref{p:HopfEinstein} implies that $[J,A_{|_Q}]=0$, that $A\xi=\lambda\xi$ on $M$ for some $\lambda\in C^\infty(M,\mathbb{R})$ and that $J$ defines a transversal K\"ahler structure on the normal bundle of the flow.
Furthermore, $\xi\left(\mathrm{tr}\left(A_{|_Q}\right)\right)=0$ because $N$ is Einstein by Proposition \ref{prop:einsteinkahler}.
The connection $1$-form $A^Q$ on $P_{\mathbb{U}_1}M\to M$ being assumed to be basic, the basic function $\theta$ must be constant by Proposition \ref{p:eqcaseKaehlersuite}.
Because $\varphi$ is the restriction of a parallel spinor field on $N$ and is assumed to be basic, Proposition \ref{p:eqcaseKaehlersuite} also implies that $\varphi$ is transversally parallel, which is also equivalent to $\theta=\lambda+\frac{1}{2}\mathrm{tr}\left(A_{|_Q}\right)$.\\
Conversely, if $\xi=-J\nu$, $[J,A_{|_Q}]=0$, $A\xi=\lambda\xi$ for some $\lambda$ and $\theta$ is constant equal to $\lambda+\frac{1}{2}\mathrm{tr}\left(A_{|_Q}\right)$, then $H_0=0$ and, as we computed in the proof of Proposition \ref{p:eqcaseKaehlersuite} -- and this computation holds true in general -- $\Omega\cdot_M\varphi=\frac{i}{2}\left(\lambda-(n+1)H\right)\varphi=-\frac{i}{2}\mathrm{tr}\left(A_{|_Q}\right)\varphi$.
Therefore,\\\\
$(n+1)^2\left(H_0^2-H^2\right)|\varphi|^2+|\Omega\cdot_M\varphi|^2+\theta^2 \vert \varphi\vert^2  -
2 i\theta \langle\varphi,\Omega\cdot_M\varphi\rangle=$
\begin{eqnarray*}
-\left(\lambda+\mathrm{tr}\left(A_{|_Q}\right)\right)^2|\varphi|^2+\frac{1}{4}\mathrm{tr}\left(A_{|_Q}\right)^2|\varphi|^2+\left(\lambda+\frac{1}{2}\mathrm{tr}\left(A_{|_Q}\right)\right)^2|\varphi|^2\\+\left(\lambda+\frac{1}{2}\mathrm{tr}\left(A_{|_Q}\right)\right)\mathrm{tr}\left(A_{|_Q}\right)|\varphi|^2=0,
\end{eqnarray*}
and thus (\ref{eq:integraline22-int}) is an equality. This concludes the proof of the proposition.
\hfill$\square$

Note that the compactness assumption on $N$ is actually no more necessary for the case where (\ref{eq:integraline22-int}) is an equality since that equality case actually holds \emph{pointwise}.
In particular, we may look for examples in non-compact K\"ahler-Einstein manifolds such as complex hyperbolic space.

In the case where $N$ has non-zero constant holomorphic sectional curvature, it already follows from a well-known result (see e.g. \cite[Theorem 2.1]{NiebergallRyan97}) that, any Hopf hypersurface must have constant principal curvature $\lambda$ -- and hence also $\mu$ -- if $A_{|_Q}=\mu {\rm Id}_Q$ and $\theta=\lambda+\frac{1}{2}{\rm tr}(A_{|_Q})$ is constant.

Theorem \ref{p:eqcaseKaehler} implies that the list of possible examples where (\ref{eq:integraline22-int}) is an e\-qua\-li\-ty cannot be long since the geometric conditions are restrictive.
Actually, the list is relatively short, at least in K\"ahler manifolds with constant holomorphic sectional curvature: by \cite[Theorem 4.1]{NiebergallRyan97} -- which summarizes \cite[Theorem 4.3]{Okumura75} (based on \cite{Ryan69}) for hypersurfaces in complex projective space and \cite{MontielRomero86} for hypersurfaces in complex hyperbolic space -- every real hypersurface in a K\"ahler manifold with constant nonvanishing holomorphic sectional curvature for which $M$ is Hopf and where $[A,J_{|_Q}]=0$ holds must be an open subset of a so-called hypersurface of type A.
Since hypersurfaces of type A in complex projective and hyperbolic spaces are completely classified (see e.g. \cite[Theorems 3.7, 3.8, 3.9, 3.13 \& 3.14]{NiebergallRyan97} and references therein), we can conclude with the following:

\begin{cor}\label{c:classifequalcaseKaehlercstholseccurv}
Let $(N^{n+2},g,J)$ be a domain with smooth nonempty connected boundary $M$ in either $\mathbb{C}\mathrm{P}^{\frac{n+2}{2}}$ or $\mathbb{C}\mathrm{H}^{\frac{n+2}{2}}$.
Assume $(M,g)$ carries a Riemannian flow given by a unit vector field $\xi$.
Let $N$ carry its canonical spin$^c$ structure, $M$ carry the spin$^c$ structure induced by the inner unit normal $\nu$ and $Q:=\xi^\perp\cap TM$ in turn carry the spin$^c$ structure induced by $\xi$.
Choose a connection $1$-form on the line bundle $K_N^{-1}{}_{|_M}=P_{\mathbb{U}_1}M\to M$ associated to the spin$^c$ structure of $Q$ with the help of a basic $\theta\in C^\infty(M,\mathbb{R})$ as above.
Assume moreover that:
\begin{enumerate}
\item the Riemannian flow $(M,g,\xi)$ is minimal,
\item there exists a section $\varphi$ of $\Sigma M$ such that $\nabla_\xi\varphi=0$ and $D_b\varphi=\frac{(n+1)H_0}{2}\varphi$ for some nonnegative basic function $H_0$ on $M$.
\end{enumerate}
Then $(n+1)^2\left(H_0^2-H^2\right)|\varphi|^2+|(\Omega\cdot_M -i \theta)\varphi|^2=0$ holds pointwise on $M$ if and only if $M$ is either
\begin{itemize}
\item a tube around a totally geodesic $\mathbb{C}\mathrm{P}^k$, where $0\leq k\leq\frac{n-2}{2}$, in case $N$ lies in $\mathbb{C}\mathrm{P}^{\frac{n+2}{2}}$, or
\item a tube around a totally geodesic $\mathbb{C}\mathrm{H}^k$, where $0\leq k\leq\frac{n}{2}$,  in case $N$ lies in $\mathbb{C}\mathrm{H}^{\frac{n+2}{2}}$.
\end{itemize}
\end{cor}

Examples include geodesic hyperspheres (case where $k=0$).
Note that all hypersurfaces of type A have two or three distinct principal curvatures.

\begin{cor} \label{coro-even}Let $(N^{n+2},g,J)$ (for $n>2$) be any compact K\"ahler manifold such that $\mathrm{Scal}_N\geq c_{n+2}|F^N|,$ where $F^N$ is the curvature form of the canonical spin$^c$ structure. Assume that the boundary $M$ of $N$ is connected of positive mean curvature $H$ and carries a minimal Riemannian flow given by a unit vector field $\xi$. Assume moreover that there exists a basic section $\varphi$ of $\Sigma M$ such that $D_b\varphi=\frac{(n+1)H_0}{2}\varphi,$ where $H_0$ is a non-negative basic function on $M$ satisfying the condition $H_0^2+\left(\frac{\sqrt{n}}{2(n+1)}|\Omega|+\frac{1}{n+1}|\theta|\right)^2\leq H^2.$ Hence, the equality case in \eqref{eq:integraline22-int} is realized and moreover $A|_Q=\mu {\rm Id}$ for some real number $\mu.$ Mainly that means the flow is transversally Einstein-K\"ahler and the manifold $M$ is $\eta$-Einstein Sasakian manifold (up to some rescaling on the metric).
\end{cor}

{\it Proof}: We estimate the two terms $|\Omega\cdot_M\varphi|^2$ and $\theta \langle\Omega\cdot_M\varphi,\varphi\rangle$ in Inequality \eqref{eq:integraline22-int}. We have
\begin{eqnarray*}
|\Omega\cdot_M\varphi|^2 &\leq &\frac{1}{4}|\sum_{i=1}^n e_i\cdot_M h(e_i)\cdot_M\varphi|^2\\
&\leq &\frac{n}{4}\sum_i |e_i \cdot_M h(e_i)\cdot_M\varphi|^2=\frac{n}{4}|h|^2|\varphi|^2=\frac{n}{2}|\Omega|^2|\varphi|^2.
\end{eqnarray*}
Recall here that we use the formula $|\Omega|^2=\frac{1}{2}|h|^2.$ For the second term, we compute
$$i\theta \langle\Omega\cdot_M\varphi,\varphi\rangle\leq |\theta||\Omega\cdot_M\varphi||\varphi|\leq |\theta|\frac{\sqrt{n}}{\sqrt{2}}|\Omega||\varphi|^2.$$
Therefore, Inequality \eqref{eq:integraline22-int} reduces to the estimate
\begin{equation}\label{eq:inequality}
0\leq\int_M\frac{1}{H}\cdot\left((n+1)^2\left(H_0^2-H^2\right)+\frac{n}{2}|\Omega|^2+\theta^2+2|\theta|\frac{\sqrt{n}}{\sqrt{2}}|\Omega|\right)|\varphi|^2dv.
\end{equation}
Since the condition $H_0^2+\left(\frac{\sqrt{n}}{2(n+1)}|\Omega|+\frac{1}{n+1}|\theta|\right)^2\leq H^2$ is fulfilled, we get the equality case in both \eqref{eq:inequality} and \eqref{eq:integraline22-int}. That means all above inequalities are sharp. In particular, one gets
$$\Omega\cdot_M\varphi=-i\frac{\sqrt{n}}{\sqrt{2}}|\Omega|\varphi.$$
Combining the last relation with the equality $\Omega\cdot_M\varphi=i(\lambda-\theta)\varphi,$ one deduces that
$$-\frac{1}{2}{\rm tr}(A|_Q)=\lambda-\theta= -\frac{\sqrt{n}}{\sqrt{2}} |\Omega|=-\frac{\sqrt{n}}{2}|h|=-\frac{\sqrt{n}}{2}|A|_Q.$$
That is, $|A|_Q=\frac{1}{\sqrt{n}}{\rm tr}(A|_Q)$ which is the equality in the Cauchy-Schwarz inequality. That yields to $A|_Q=\mu {\rm Id}$ for some function $\mu.$ With the help of Proposition \ref{p:hypersisopar}, one deduces that both $\mu$ and $\lambda$ are constant. In view of Proposition \ref{p:bordSasaki}, the manifold $M$ is a Sasakian $\eta$-Einstein manifold from the fact that the transverse metric remains invariant.
\hfill$\square$

\subsection{The odd-dimensional case}\label{subsec:odddimcase}

In this section, we look at the case where the normal bundle has odd rank.
As mentioned above, we will be interested in the family of manifolds $N$ that are domains in Riemannian products of K\"ahler-Einstein manifolds $N_1$ with $\R$ or $\mathbb{S}^1.$

Let $M$ be the boundary of any domain $N$ in $N_1\times\R$ or $N_1\times \mathbb{S}^1$ and carry a Riemannian flow.
Since the rank of the normal bundle $Q$ is odd, we have the identification
$$\Sigma Q\oplus \Sigma Q\simeq \Sigma M\simeq \textbf{S},$$
where in the first isomorphism, we use the following identifications for the Clifford multiplications
$$(Z\cdot_Q\oplus -Z\cdot_Q)\Upsilon=\xi\cdot_M Z\cdot_M\Upsilon,$$
for any $Z\in\Gamma(Q)$ and $\Upsilon\in \textbf{S}=\Sigma N|_M.$
For the second isomorphism, we have $X\cdot_M\Upsilon=X\cdot\nu\cdot\Upsilon.$
Now the action of $i\nu$ on $\textbf{S}$ is determined by the action of the complex volume form $\omega$ of $\Sigma M$, that is for any spinor $\Upsilon$ on $\textbf{S},$ we have $i\nu\cdot\Upsilon=\omega\cdot\Upsilon=\bar{\Upsilon}$, where $\bar{\Upsilon}=\Upsilon_+ -\Upsilon_{-}$ with $\Upsilon_\pm$ are eigensections of $\omega$ corresponding to the eigenvalues $\pm 1.$
Thus, from the definition of the projections $P_\pm$, we deduce that $P_\pm\Upsilon=\Upsilon_\pm$.


{\it Proof of Theorem \ref{pro:estimate-int} for $n$ odd}: Let us define the spinor field $\Upsilon=\varphi+\xi\cdot_M\varphi$
where $\varphi$ is considered as a section in $\Sigma Q\simeq
\Sigma^+M,$ i.e. $P_+\varphi=\varphi.$
Using  \eqref{eq:relationdirac}, we have
\begin{eqnarray} \label{eq:relations}
{\textbf D}_{\textbf{S}}\Upsilon=D_M\Upsilon&=&D_M\varphi+D_M(\xi\cdot_M\varphi)\nonumber\\
&=& -\frac{n+1}{2}H_0\xi\cdot_M\varphi-\frac{1}{2}\xi\cdot_M\Omega\cdot_M\varphi +{\frac i2 \theta \xi\cdot_M\varphi } \nonumber \\
&& +\xi\cdot_M D_b(\xi\cdot_M\varphi)+\frac{1}{2}\Omega\cdot_M\varphi {- \frac i2 \theta \varphi}\nonumber\\
&\bui{=}{\eqref{eq:commute}}&
-\frac{n+1}{2}H_0\xi\cdot_M\varphi-\frac{1}{2}\xi\cdot_M\Omega\cdot_M\varphi
+ \frac i2 \theta \xi\cdot_M\varphi \nonumber \\
&& -\frac{n+1}{2}H_0\varphi+\frac{1}{2}\Omega\cdot_M\varphi
{- \frac i2 \theta \varphi}.
\end{eqnarray}

It is easy to check from $\varphi\in \Gamma(\Sigma^+M)$ that the identities $\langle \xi\cdot_M\Omega\cdot_M\varphi,
\varphi \rangle=0$ and $\langle\xi\cdot_M\varphi, \varphi\rangle=0$ hold.
Hence by taking the norm of the spinor field $|{\textbf D}_{\textbf{S}}\Upsilon|^2,$ we find that it is equal to
\begin{eqnarray*}
|{\textbf D}_{\textbf{S}}\Upsilon|^2 =
\frac{(n+1)^2}{2}H_0^2|\varphi|^2+\frac{1}{2}|\Omega\cdot_M\varphi|^2
+\frac{\theta^2}{2} \vert \varphi\vert^2-
i\theta \langle\varphi, \Omega\cdot_M\varphi\rangle.
\end{eqnarray*}
Inequality \eqref{hij-mon} applied to the spinor field $\Upsilon$ finishes the proof by plugging the last equality  and
using the fact that $|\Upsilon|^2=2|\varphi|^2.$
\hfill$\square$

The proof of Theorem \ref{pro:caslimiteimpair} is technical and will be splitted into several lemmas (Lemmas \ref{lemma:equ1}, \ref{lem:cara1}, \ref{lem:carac1}, \ref{lem:carac2}, \ref{lem:caract12}, \ref{lem:caracte3}). First, we have

\begin{lemma} \label{lemma:equ1} If the equality case is realized in \eqref{eq:integraline22-int}, then
 $$-h(X)\cdot_M\varphi+(\frac{1}{2b}-\frac{b}{2})A(X)\cdot_M\varphi=\frac{1}{b}g(A(X),\xi)\xi\cdot_M\varphi,$$
for all $X\in \Gamma(TM)$ and for some $b\in\mathbb{C}.$
In particular, we get $(\frac{1}{2b}-\frac{b}{2})^2|A(\xi)|^2=-g(A(\xi),\xi)^2.$
\end{lemma}

\noindent{\it Proof}: Assume that the equality is realized in \eqref{eq:integraline22-int} and recall that $\Upsilon=\varphi+\xi\cdot_M\varphi,$ then by Proposition \ref{pro:hijazimontiel} there exists two parallel spinors $\psi$ and $\vartheta$ in $\Sigma N$ such that $P_+\Upsilon=\varphi=\psi_+$ and $P_-\Upsilon=\xi\cdot_M\varphi=\vartheta_-.$ Since the dimension of the space of parallel spinors on $N_1\times \mathbb{S}^1$ is one (those are being identified with parallel spinors on $N_1$ corresponding to the canonical spin$^c$ structure), we deduce that $\vartheta=b\psi$ for some $b\in \mathbb{C}.$ Now differentiating the equation $\xi\cdot_M\varphi=b\psi_-$ in the direction of any vector field $X\in \Gamma(TM),$ we get that
$$
h(X)\cdot_M\varphi-\frac{1}{2}\xi\cdot_M A(X)\cdot_M\psi_-=-\frac{b}{2}A(X)\cdot_M\varphi.
$$
Here, we used the fact that $\nabla^M_X\psi_\pm=-\frac{1}{2}A(X)\cdot_M\psi_\mp,$ since $\psi$ is parallel on $N.$ Replacing now $\psi_-$ by $\frac{1}{b}\xi\cdot_M\varphi$, the above equation reduces to
\begin{equation}\label{eq:limitcaseodd}
-h(X)\cdot_M\varphi+(\frac{1}{2b}-\frac{b}{2})A(X)\cdot_M\varphi=\frac{1}{b}g(A(X),\xi)\xi\cdot_M\varphi.
\end{equation}
To prove the second part of the lemma, we take $X=\xi$ in \eqref{eq:limitcaseodd} and use the fact that the mean curvature $\kappa=\nabla_\xi^M\xi$ vanishes to deduce that
\begin{eqnarray} \label{eq:limitcaseconsequence}
(\frac{1}{2b}-\frac{b}{2})A(\xi)\cdot_M\varphi=\frac{1}{b}g(A(\xi),\xi)\xi\cdot_M\varphi.
\end{eqnarray}
By taking the Clifford multiplication of \eqref{eq:limitcaseconsequence} by $\xi$ and applying the rule
$\xi\cdot_M A(\xi)\cdot_M=-A(\xi)\cdot_M\xi\cdot_M-2g(A(\xi),\xi),$ we get that
\begin{eqnarray} \label{eq:limitcaseconsequencebis}
(\frac{1}{2b}-\frac{b}{2})A(\xi)\cdot_M\xi\cdot_M\varphi=bg(A(\xi),\xi)\varphi.
\end{eqnarray}
Now if $g(A(\xi),\xi)=0,$ then either $(\frac{1}{2b}-\frac{b}{2})=0$ or $A(\xi)=0$ from which the relation in the lemma is proved. If $g(A(\xi),\xi)$ is different from zero, we use again Equation \eqref{eq:limitcaseconsequence} divided by the term $\frac{1}{b}g(A(\xi),\xi)$ to replace $\xi\cdot_M\varphi$ by its value in \eqref{eq:limitcaseconsequencebis}. The result then follows.
\hfill$\square$

\begin{lemma} \label{lem:cara1} If the equality case is realized in \eqref{eq:integraline22-int} and $M$ is connected, then $g(\nu,\partial_t)=0, |b|=1.$
Moreover, there exists a smooth compact domain $\Delta$ with boundary $M_1$ in $N_1$ such that $N=\Delta\times\mathbb{S}^1$, in particular $M=M_1\times\mathbb{S}^1$.
\end{lemma}

\noindent{\it Proof}: As in the previous lemma, we know that $\varphi=\psi_+$ and $\xi\cdot_M\varphi=b\psi_-,$ where $\psi$ is a parallel spinor of norm assumed to be equal to $1$. Therefore, we deduce that $|\varphi|^2=|\psi_+|^2=\frac{|b|^2}{|b|^2+1}$ and that $|\psi_-|^2=\frac{1}{|b|^2+1}.$ Hence, we get that
$$\langle i\nu\cdot\psi,\psi\rangle=|\psi_+|^2-|\psi_-|^2=\frac{|b|^2-1}{|b|^2+1}.$$

\noindent Now, for every $X\in TN=TN_1\oplus \R\partial_t$, we may split $X=X^T+g(X,\partial_t)\partial_t$, where $X^T$ is pointwise tangent to the $N_1$-factor. In particular $X\cdot\psi=X^T\cdot\psi+g(X,\partial_t)\partial_t\cdot\psi$; but because of $\psi\in\Sigma_0 N_1$ pointwise, we have $X^T\cdot\psi=-iJ(X^T)\cdot\psi$, so that, using also $i\partial_t\cdot\psi=\psi$ (because of $\psi\in\Gamma(\Sigma_+N_1)$ pointwise), we get
\begin{equation}\label{eq:identity}
X\cdot\psi=-iJ(X^T)\cdot\psi-ig(X,\partial_t)\psi.
\end{equation}

\noindent As a consequence, $\langle iX\cdot\psi,\psi\rangle=\langle J(X^T)\cdot\psi,\psi\rangle+g(X,\partial_t)|\psi|^2$; but since both $\langle iX\cdot\psi,\psi\rangle$ and $g(X,\partial_t)|\psi|^2$ are real whereas $\langle J(X^T)\cdot\psi,\psi\rangle$ is purely imaginary, we deduce that in fact $\langle J(X^T)\cdot\psi,\psi\rangle=0$ and $\langle iX\cdot\psi,\psi\rangle=g(X,\partial_t)|\psi|^2=g(X,\partial_t)$.
This implies first that $g(\nu,\partial_t)=g(\nu,\partial_t)|\psi|^2=\langle i\nu\cdot\psi,\psi\rangle=\frac{|b|^2-1}{|b|^2+1}$, in particular $g(\nu,\partial_t)$ is constant on $M$.
Now if $M$ is connected, then we may apply the divergence theorem and obtain for the parallel vector field $\partial_t$ on $N$
\[0=\int_N\delta^N(\partial_t) d\mu_g^N=\pm\int_M g(\nu,\partial_t) d\mu_g^M,\]
from which $g(\nu,\partial_t)=0$ follows.
In particular, $|b|=1$.\\

It remains to show the existence of a domain $\Delta$ of $N_1$ such that $N=\Delta\times \mathbb{S}^1$ (and hence $M=M_1\times\mathbb{S}^1$, where $M_1:=\partial\Delta$).
For this, we show that, for any $t\in\R$, the flow $\phi_t$ of $\partial_t$ preserves $N$, that is, that $\phi_t(N)=N$.
First consider the case where $N\subset N_1\times\R$.
We may assume that $0\in t(N)$ and identify $N_1$ with $N_1\times\{0\}\subset N_1\times\R$, so that $N_1$ becomes the preimage of the regular value $0$ for the function $t$ on $N_1\times\R$.
Since $\partial_t$ is parallel on $N_1\times\R$ and $N_1$ is assumed to be complete, so is $N_1\times\R$ and the flow $\phi$ of $\partial_t$ is defined on $N_1\times\R$; actually, $\phi$ is the identity map on $N_1\times\R$.
Moreover, because the restriction of $\partial_t$ onto $M$ is tangent to $M$, the flow $\phi_s$ preserves $M$ for all $s\in\R$, that is, $\phi_s(M)=M$.
Let $x\in N$, then either $x\in \partial N=M$ and then $\phi_t(x)\in M\subset N$ as we have just seen; or $x\in\bui{N}{\circ}$, but then $\phi_t(x)$ can lie neither on $M$ (otherwise $x=\phi_{-t}(\phi_t(x))\in M$) nor outside $N$ (otherwise the integral curve $s\mapsto\phi_s(x)$ linking $x$ with $\phi_t(x)$ must cross $M$ and thus lie in $M$ for all time), therefore $\phi_t(x)\in\bui{N}{\circ}\subset N$.
On both cases, $\phi_t(x)\in N$.
Therefore $\phi_t(N)\subset N$; changing $t$ into $-t$ gives $N\subset\phi_t(N)$ and hence $N=\phi_t(N)$.\\
Now because $\phi_t$ preserves $N$ for all $t\in\R$, we may set $\Delta:=t^{-1}(\{0\})\cap N\subset N$, which is a smooth domain with boundary $M_1=t^{-1}(\{0\})\cap M$ in $N_1\times \R$ (it is smooth up to the boundary because of $g(\nu,\partial_t)=0$).
Because $\partial_t$ is parallel on $N$, the flow $\phi$ induces an isometry $\Delta\times\R\to N$, in particular $\Delta$ must be connected.
Since $\phi$ is the identity on $N_1\times\R$, we have actually shown that $N=\Delta\times\R$.
Note that this case cannot happen here since $N$ is assumed to be compact.\\
In case where $N\subset N_1\times \mathbb{S}^1$, we may lift $N$ to a smooth manifold with boundary $\overline{N}$ in $N_1\times\R$ via the covering map $N_1\times\R\to N_1\times \mathbb{S}^1$ -- which is not the universal cover of $N_1\times \mathbb{S}^1$, unless $N_1$ itself is simply-connected.
Since that covering map preserves $\partial_t$, that parallel vector field is tangent to $\overline{M}=\partial\overline{N}$ and therefore $\overline{N}=\Delta\times\R$ for some smooth domain $\Delta$ in $N_1$ by the above argument.
Now because the $\mathbb{Z}$-group action on $N_1\times\R$ underlying the covering map is trivial on $N_1$ (it only acts on the $\R$-factor), we can conclude that $N=\Delta\times \mathbb{S}^1$.
This concludes the proof.
\hfill$\square$


The computation in the sequel will be devoted to show that in the equality case of \eqref{eq:integraline22-int}, the vector field $\xi$ defining the flow will be equal to $\pm\partial_t.$ The main idea is to show that the $TM_1$-component of $\xi$ corresponding to the decomposition $TM=TM_1\oplus \mathbb{R}\partial_t$ (according to Lemma \ref{lem:cara1}) defines a Riemannian flow on $M_1$ and the solution $\varphi$ of the basic Dirac equation defines as well a solution of the basic Dirac equation corresponding to that flow (which has even-dimensional normal bundle). It turns out that such a solution realizes the equality case of the integral inequality established in the previous section. We begin with a remark of algebraic nature:

\begin{lemma}\label{l:actionpsi+vanish}
Let $Z,W\in TN$ be such that $(Z+iW)\cdot\psi_+=0$, where $\psi$ comes from a nonzero constant section of $\Sigma_0N_1$.
Then
\[Z^T+J(W^T)-g(Z,\partial_t)\nu-g(W,\partial_t)J(\nu)=0,\]
where $(\cdot)^T$ is the orthogonal projection onto $TN_1$.
\end{lemma}

{\it Proof:}
Because of $\psi\in\Sigma_0N_1$ and hence $i\partial_t\cdot\psi=\psi$, we have $X\cdot\psi=X^T\cdot\psi+g(X,\partial_t)\partial_t\cdot\psi=X^T\cdot\psi-ig(X,\partial_t)\psi$ for every $X\in TN$.
As a consequence, we can write
\begin{eqnarray*}
Z\cdot\psi_+&=&\frac{1}{2}\left(Z\cdot\psi+Z\cdot i\nu\cdot\psi\right)\\
&=&\frac{1}{2}\left(Z^T\cdot\psi-ig(Z,\partial_t)\psi+iZ\cdot \nu^T\cdot\psi+\underbrace{g(\nu,\partial_t)}_{0}Z\cdot\psi\right)\\
&=&\frac{1}{2}\left(Z^T\cdot\psi-ig(Z,\partial_t)\psi+iZ^T\cdot \nu^T\cdot\psi+ig(Z,\partial_t)\partial_t\cdot\nu^T\cdot\psi\right)\\
&=&\frac{1}{2}\left(Z^T\cdot\psi-ig(Z,\partial_t)\psi+iZ^T\cdot \nu^T\cdot\psi-g(Z,\partial_t)\nu^T\cdot\psi\right),
\end{eqnarray*}
which yields
\begin{eqnarray*}
\left(Z+iW\right)\cdot\psi_+&=&\frac{1}{2}\left(Z^T\cdot\psi-ig(Z,\partial_t)\psi+iZ^T\cdot \nu^T\cdot\psi-g(Z,\partial_t)\nu^T\cdot\psi\right)\\
&&+\frac{i}{2}\left(W^T\cdot\psi-ig(W,\partial_t)\psi+iW^T\cdot \nu^T\cdot\psi-g(W,\partial_t)\nu^T\cdot\psi\right)\\
&=&\frac{1}{2}\Big((Z^T+iW^T)\cdot\psi-ig(Z+iW,\partial_t)\psi\\
&&\phantom{\frac{1}{2}\Big(}+i(Z^T+iW^T)\cdot\nu^T\cdot\psi-g(Z+iW,\partial_t)\nu^T\cdot\psi\Big),
\end{eqnarray*}
where $g$ has been extended as a complex bilinear form on $TN\otimes\mathbb{C}$.
Now we can make use of identity \eqref{eq:identity}, that uses specifically $\psi\in\Sigma_0N_1$: for every $X\in TN$, $X\cdot\psi=-iJ(X^T)\cdot\psi-ig(X,\partial_t)\psi$.
Since $(Z+iW)\cdot\psi_+=0$ if and only if $\partial_t\cdot(Z+iW)\cdot\psi_+=0$, we obtain that $(Z+iW)\cdot\psi_+=0$ is equivalent to
\begin{eqnarray*}
0&=&\partial_t\cdot(Z^T+iW^T)\cdot\psi-ig(Z+iW,\partial_t)\partial_t\cdot\psi\\
&&+i\partial_t\cdot\underbrace{(Z^T+iW^T)\cdot_{N_1}\nu^T\cdot_{N_1}\psi}_{\in\Sigma_+N_1}-g(Z+iW,\partial_t)\partial_t\cdot\nu^T\cdot\psi\\
&=&-(Z^T+iW^T)\cdot_{N_1}\psi-g(Z+iW,\partial_t)\psi+(Z^T+iW^T)\cdot_{N_1}\nu^T\cdot_{N_1}\psi\\
&&+g(Z+iW,\partial_t)\nu^T\cdot_{N_1}\psi\\
&=&-(Z^T+J(W^T))\cdot_{N_1}\psi-g(Z+iW,\partial_t)\psi+(Z^T+iW^T)\cdot_{N_1}\nu^T\cdot_{N_1}\psi\\
&&+g(Z,\partial_t)\nu^T\cdot_{N_1}\psi+g(W,\partial_t)J(\nu^T)\cdot_{N_1}\psi.
\end{eqnarray*}
Identifying the degree-$1$-terms on both sides and using $\psi\neq0$, we obtain $0=-(Z^T+J(W^T))+g(Z,\partial_t)\nu^T+g(W,\partial_t)J(\nu^T)$, which by $\nu=\nu^T$ is the result.
\hfill$\square$

\begin{lemma}\label{l:xiTisproporttoJnu}
If $g(\nu,\partial_t)=0$, then $(\xi+b\partial_t)\cdot\varphi=0$, from which follow that $g(\xi,\partial_t)=-\Re(b)$ and
\begin{equation}\label{eq:xiTJnu}
\xi^T=\Im m(b)J\nu.
\end{equation}
In particular, $hX=-\Im m (b)J(AX)^{TM}$ for all $X\in TM$.
\end{lemma}

{\it Proof:} Since $\nu$ is orthogonal to $\partial_t$, one can easily check that $i\partial_t\cdot P_+\psi$ belongs to $\Sigma_-M$ (here $\Sigma_-M$ denotes the eigenspace of the action of $i\nu$ corresponding to the eigenvalue $-1$).
Hence, by using the fact that $i\partial_t\cdot\psi=\psi$ and the decomposition $\psi=P_+\psi+P_-\psi,$ where $P_+\psi=\varphi$ and $P_-\psi=\frac{1}{b}\xi\cdot_M\varphi,$  we deduce that
\begin{equation}\label{eq:relationxipartial}
\partial_t\cdot\varphi=-\frac{1}{b}\xi\cdot\varphi,
\end{equation}
which gives the first identity.
Setting $Z:=\xi+\Re(b)\partial_t$ and $W:=\Im m(b)\partial_t$, this identity becomes $(Z+iW)\cdot\psi_+=0$, so that $Z^T+J(W^T)-g(Z,\partial_t)\nu-g(W,\partial_t)J(\nu)=0$ by Lemma \ref{l:actionpsi+vanish}.
Because of $Z^T=\xi^T$ and $W^T=0$, this identity is equivalent to $\xi^T-(g(\xi,\partial_t)+\Re(b))\nu-\Im m(b)J\nu=0$, which itself is equivalent to $g(\xi,\partial_t)=-\Re(b)$ (thing we can anyway read off the real part of the inner product of $(\xi+b\partial_t)\cdot\varphi=0$ with $\varphi$) and $\xi^T=\Im m(b)J\nu$, which is \eqref{eq:xiTJnu}.
As an immediate consequence, for every $X\in TM$, $hX=\nabla_X^M\xi=\nabla_X^M\xi^T=\Im m(b)\nabla_X^M(J\nu)=-\Im m(b)J(AX)^{TM}$.
\hfill$\square$


\begin{lemma}\label{lem:carac1}
Assume that the equality in \eqref{eq:integraline22-int} is realized, then by writing $b=e^{i\beta}$ for some $\beta\in \mathbb{R},$ there exists an $\varepsilon\in\{\pm1\}$ such that 
$$|A(\xi)|{\rm sin}\beta=\varepsilon g(A(\xi),\xi) \ \ \text{and}\ \ \  \theta=-\frac{1}{2}((n+1)H {\rm sin}\beta+\varepsilon|A(\xi)|).$$
Moreover, ${\rm cos}\beta=-\frac{H_0}{H}$.
\end{lemma}

\noindent{\it Proof}: Because of $|b|=1$, we have $\Im m(b)^2|A\xi|^2=g(A\xi,\xi)^2$ by Lemma \ref{lemma:equ1}, so that there is an $\varepsilon\in\{\pm1\}$ such that $|A(\xi)|{\rm sin}\beta=\varepsilon g(A(\xi),\xi)$.

By tracing \eqref{eq:limitcaseodd} over an orthonormal frame on $\Gamma(TM)$, we get that
$$-2\Omega\cdot_M\varphi+i(n+1)H {\rm sin}\beta\varphi=\frac{1}{b}A(\xi)\cdot_M\xi\cdot_M\varphi.$$
But by \eqref{eq:limitcaseconsequencebis}, $\frac{1}{b}A\xi\cdot_M\xi\cdot_M\varphi=i\varepsilon|A\xi|\varphi$ (this holds true whether $A\xi$ vanishes or not), so that
\begin{equation}\label{eq:actionOmegaodd}
\Omega\cdot_M\varphi=\frac{i}{2}\left((n+1)H\Im m(b)-\varepsilon|A\xi|\right)\varphi.
\end{equation}
In order to compute $\theta$, we use the first equation in \eqref{eq:oneillspin} and the fact that $\varphi$ is a basic spinor.
In fact,
\begin{eqnarray*}
\frac{i\theta}{2}\varphi&=&\nabla_\xi^M\varphi-\frac{1}{2}\Omega\cdot_M\varphi\\
&\bui{=}{\eqref{eq:actionOmegaodd}}&-\frac{1}{2}A(\xi)\cdot_M\psi_--\frac{i}{4}\left((n+1)H\Im m(b)-\varepsilon|A\xi|\right)\varphi\\
&=&-\frac{1}{2b}A(\xi)\cdot_M\xi\cdot_M\varphi-\frac{i}{4}\left((n+1)H\Im m(b)-\varepsilon|A\xi|\right)\varphi\\
&\bui{=}{\eqref{eq:limitcaseconsequencebis}}&-\frac{i\varepsilon}{2}|A(\xi)|\varphi-\frac{i}{4}\left((n+1)H\Im m(b)-\varepsilon|A\xi|\right)\varphi\\
&=&-\frac{i}{4}\left((n+1)H\Im m(b)+\varepsilon|A\xi|\right)\varphi,
\end{eqnarray*}
which gives the expression of $\theta$.
To compute the real part of $b,$ we use the second relation in \eqref{eq:relationdirac} between the extrinsic Dirac operator of $M$ and the basic Dirac operator and the fact that $D_b\varphi=\frac{(n+1)H_0}{2}\varphi$.
Indeed,
$$\textbf{D}_\textbf{S}\varphi=D_M\varphi=-\frac{n+1}{2}H_0\xi\cdot_M\varphi-\frac{1}{2}\xi\cdot_M\Omega\cdot_M\varphi+\frac{i\theta}{2}\xi\cdot_M\varphi.$$
But recall that $\varphi=P_+\psi$ on $M$ where $\psi$ is parallel on $N.$ Therefore, the extrinsic Dirac operator applied to $\varphi$ is equal to
$$\textbf{D}_\textbf{S}\varphi=\textbf{D}_\textbf{S}(P_+\psi)=P_-(\textbf{D}_\textbf{S}\psi)=\frac{n+1}{2}HP_-\psi=\frac{n+1}{2b}H\xi\cdot_M\varphi.$$
Hence by comparing the above two equalities, we deduce that $(n+1)(H_0+\Re(b)H)=0$, from which $\Re(b)=-\frac{H_0}{H}$ follows.
\hfill$\square$

\begin{lemma}\label{lem:carac2}
If \eqref{eq:xiTJnu} holds, then $\theta$ is constant on $M$.
\end{lemma}

{\it Proof:} Using \eqref{eq:AQbasique} along any vector field: for any $X\in \Gamma(TM)$
$$F^M(\xi,X)=F^N(\xi,X)=-i{\rm Ric}^{N_1}(J\xi^T,X)=-i\frac{{\rm Scal}^{N_1}}{n+1}g(J(\xi^T),X),$$
so that $\xi\lrcorner F^M=-i\frac{{\rm Scal}^{N_1}}{n+1}J(\xi^T)^{TM}=0$ by \eqref{eq:xiTJnu}.
\hfill$\square$

It is important to notice here that, computing the integrand on the r.h.s. of \eqref{eq:integraline22-int} and using what we already know about $\theta$, $\Omega$, $H_0$ and $H$, we find
\begin{eqnarray*}
&\ & (n+1)(H_0^2-H^2)|\varphi|^2+\left|(\Omega-i\theta)\cdot_M\varphi\right|^2 \\ 
&=&  (n+1)^2(\cos(\beta)^2-1)H^2|\varphi|^2 +(n+1)^2H^2\sin^2(\beta)|\varphi|^2\\
&=&0,
\end{eqnarray*}
in particular  \eqref{eq:integraline22-int} is an equality.
However, this does not suffice to conclude.
We shall indeed show that, if $\Im m(b)\neq0$, then we are led to a contradiction.\\

From now on we assume that $\Im m(b)\neq0$, that is, that $b\neq\pm1$.
By Lemma \ref{l:xiTisproporttoJnu}, if $\Im m(b)\neq0$, then $0=h\xi=-\Im m(b)J(A\xi)^{TM}$ implies that $A\xi=A\xi^T$ is proportional to $\xi^T$.
More precisely, $A(\xi)=\varepsilon\frac{|A(\xi)|}{{\rm sin}\beta}\xi^T$.
Mainly, that means the vector field $\xi_1:=\xi^T$ is an eigenvector of $A$ corresponding to the eigenvalue $\lambda_1:=\varepsilon\frac{|A(\xi)|}{{\rm sin}\beta}.$
Moreover, the vector field $\xi_1$ is of constant norm equal to $|{\rm sin}\beta|$ and defines a minimal Riemannian flow on the manifold $M_1$, isometric to the product, with O'Neill tensor
$$h_1:=\nabla^{M_1}(\frac{\xi_1}{|\xi_1|})=\frac{1}{|{\rm sin}\beta|}\nabla^M\xi.$$

\noindent The manifold $M_1$ is clearly spin$^c$ with a connection form $A^{M_1}=A^M|_{M_1}.$ Hence as mentioned in Section \ref{sec:prel}, the normal bundle $Q_1$ carries also a spin$^c$ structure with the same line bundle as for $M_1.$ Now, we choose a connection $1$-form on $Q_1$ as
$$A^{Q_1}:=A^{M_1}-i\theta_1\frac{\xi_1}{|\xi_1|},$$
where $\theta_1:=\frac{\theta}{|\xi_1|}.$ The relation \eqref{eq:AQbasique} is clearly satisfied on $M_1,$ since
$$(\frac{\xi_1}{|\xi_1|})\lrcorner F^{M_1}=(\frac{\xi}{|\xi_1|})\lrcorner F^{M}=-id(\frac{\theta}{|\xi_1|})=-id\theta_1.$$
Also, one can check by choosing ${\rm sin}\beta<0$ that $\frac{\xi_1}{|\xi_1|}=-J\nu$ and from the expression of $\theta$ in Lemma \ref{lem:carac1}, that
$$\lambda_1-2\theta_1=\varepsilon\frac{|A(\xi)|}{{\rm sin}\beta}-\frac{1}{{\rm sin}\beta}(\varepsilon|A\xi|+(n+1)H{\rm sin}\beta)=-(n+1)H=-nH_1,$$
where $H_1$ is the mean curvature of $M_1$ into $N_1.$ Taking into account those observations, we get the following lemma

\begin{lemma} \label{lem:caract12} Assume that the equality in \eqref{eq:integraline22-int} is realized and $b\neq \pm 1$, then the equality case of the Inequality \eqref{eq:integraline22-int} on $M_1$ is realized.
\end{lemma}

\noindent {\it Proof}: According to Theorem \ref{p:eqcaseKaehler} in Subsection \ref{subsec:evendimcase}, it is sufficient to prove the existence of a solution of the basic Dirac equation on $M_1.$ Let us denote by $\varphi_1(x):=\varphi(x,1)$ for $x\in M_1,$ where $\varphi$ is the solution of the basic Dirac equation on $M$ realizing the equality case in \eqref{eq:integraline22-int}. In the following, we aim to show that $\varphi_1$ is basic with respect to the flow on $M_1$ and is a solution of the basic Dirac equation. For this, we use the first equation in \eqref{eq:oneillspin} to compute
\begin{eqnarray*}
\nabla_{(\frac{\xi_1}{|\xi_1})}\varphi_1&=&\nabla_{(\frac{\xi_1}{|\xi_1})}^{M_1}\varphi_1-\frac{1}{2}\Omega^1\cdot_{M_1}\varphi_1-\frac{i\theta_1}{2}\varphi_1\\
&=&\frac{1}{|\xi_1|}(\nabla_\xi^M\varphi_1-\frac{1}{2}\Omega\cdot_{M}\varphi_1-\frac{i\theta}{2}\varphi_1)=0.
\end{eqnarray*}
Here, we used the fact that $\varphi_1$ is basic on $M$ and constant along $\partial_t.$ We also mention that for an orthonormal frame $\{e_i^1\}_{i=1,\cdots,n-1}$ of $Q_1\subset Q,$ we have
\begin{eqnarray*}
\Omega^1\cdot_{M_1}\varphi_1&=&\frac{1}{2}\sum_{i=1}^{n-1}e_i^1\cdot_{M_1}h_1(e_i^1)\cdot_{M_1}\varphi_1\\
&=& \frac{1}{2|\xi_1|}\sum_{i=1}^ne_i\cdot_{M}h(e_i)\cdot_{M}\varphi_1\\
&=&\frac{1}{|\xi_1|}\Omega\cdot_{M}\varphi_1
\end{eqnarray*}
where $\{e_i\}_{i=1,\cdots,n}$ is an orthonormal basis of $Q$ defined by $\{Z,e_i^1\}_{i=1,\cdots,n-1}$ and $Z$ is a linear combination of $\xi_1$ and $\partial_t,$ which gives $h(Z)=0.$

\noindent The last part would be to compute the basic Dirac operator $D_b^1$ to $\varphi_1.$ Using Equation \eqref{eq:relationdirac}, we write
\begin{eqnarray*}
D_b^1\varphi_1&=&D_{M_1}\varphi_1+\frac{1}{2|\xi_1|}\xi_1\cdot_{M_1}\Omega^1\cdot_{M_1}\varphi_1-\frac{i\theta_1}{2|\xi_1|}\xi_1\cdot_{M_1}\varphi_1\\
&=&-\partial_t\cdot_M D_M\varphi_1+\frac{1}{2|\xi_1|}\xi_1\cdot_{M_1}\Omega^1\cdot_{M_1}\varphi_1-\frac{i\theta_1}{2|\xi_1|}\xi_1\cdot_{M_1}\varphi_1 \\
&\bui{=}{\eqref{eq:relationdirac}}&\frac{n+1}{2}H_0\partial_t\cdot\xi\cdot\varphi_1+\frac{1}{2}\partial_t\cdot\xi\cdot\Omega\cdot\varphi_1-\frac{i\theta}{2}\partial_t\cdot\xi\cdot\varphi_1\\
&&+\frac{1}{2|\xi_1|^2}\xi_1\cdot\partial_t\cdot\Omega\cdot\varphi_1-\frac{i\theta_1}{2|\xi_1|}\xi_1\cdot\partial_t\cdot\varphi_1\\
&\bui{=}{\eqref{eq:relationxipartial}}&\frac{n+1}{2}H_0 b\varphi_1+\frac{b}{2}\Omega\cdot\varphi_1-\frac{ib\theta}{2}\varphi_1+\frac{1}{2{\rm sin}^2\beta}(\xi+{\rm cos}\beta\partial_t)\cdot\partial_t\cdot\Omega\cdot\varphi_1\\&&-\frac{i\theta}{2{\rm sin}^2\beta}(\xi+{\rm cos}\beta\partial_t)\cdot\partial_t\cdot\varphi_1\\
&=& \left(\frac{n+1}{2}H_0 b-\frac{ib\theta}{2}-\frac{i\theta}{2b{\rm sin}^2\beta}+\frac{i\theta{\rm cos}\beta}{2{\rm sin}^2\beta}\right)\varphi_1\\
&&+\left(\frac{b}{2}+\frac{1}{2b{\rm sin}^2\beta}-\frac{{\rm cos}\beta}{2{\rm sin}^2\beta}\right)\Omega\cdot\varphi_1\\
\end{eqnarray*}
Then, using \eqref{eq:actionOmegaodd}, we get
\begin{eqnarray*}
D_b^1\varphi_1&=&\left(\frac{n+1}{2}H_0 b-\frac{ib\theta}{2}-\frac{i\theta}{2b{\rm sin}^2\beta}+\frac{i\theta{\rm cos}\beta}{2{\rm sin}^2\beta}\right)\varphi_1\\
&&+\left(\frac{b}{2}+\frac{1}{2b{\rm sin}^2\beta}-\frac{{\rm cos}\beta}{2{\rm sin}^2\beta}\right)\frac{i}{2}((n+1)H{\rm sin}\beta-\varepsilon|A(\xi)|)\varphi_1.
\end{eqnarray*}
By writing $b={\rm cos}\beta+i{\rm sin}\beta$ and using the fact that $-\varepsilon|A(\xi)|=2\theta+(n+1)H{\rm sin}\beta$ and that ${\rm cos}\beta=-\frac{H_0}{H}$ from Lemma \ref{lem:carac1}, we find after a straitforward computation that $D_b^1\varphi_1=0.$
\hfill$\square$

Now, we deduce with the following lemma:

\begin{lemma} \label{lem:caracte3} Assume that the equality in \eqref{eq:integraline22-int} is realized, then $b=\pm 1$.
\end{lemma}

\noindent {\it Proof}: Assume that $b\neq \pm 1,$ then from Lemma \ref{lem:caract12} the manifold $M_1$ is a limiting manifold for the even case. As a consequence the spinor field $\varphi_1$ is the restriction of a parallel spinor on $N_1$ and that $i\frac{\xi_1}{|\xi_1|}\cdot_{M_1}\varphi_1=\varphi_1.$ Then, we write
$i(\xi+{\rm cos}\beta\partial_t)\cdot_M\partial_t\cdot_M\varphi_1=-{\rm sin}\beta\varphi_1.$ That gives the following identity $\xi\cdot\partial_t\cdot\varphi_1=b\varphi_1.$ Now combining the last relation with Equation $(\xi+b\partial_t)\cdot\varphi=0$ gives that $b^2=1.$
\hfill$\square$

\noindent{\it Proof of Theorem \ref{pro:caslimiteimpair}}: As $b=\pm 1$ from Lemma \ref{lem:caracte3}, we deduce from Equation $(\xi+b\partial_t)\cdot\varphi=0$ that $\xi=\pm\partial_t.$ In particular, the first equation in \eqref{eq:oneillspin} implies that $\theta=0$ and the computation in Lemma \ref{lem:carac1} gives that $H_0=H.$

For the converse, assume that $N$ is isometric to $\Delta\times \mathbb{S}^1$ where $\Delta$ is a K\"ahler-Einstein manifold with boundary $M_1$  and let $\xi=\partial_t$ be the parallel vector field that defines the Riemannian flow on the boundary $M=M_1\times \mathbb{S}^1$ (that is, $h=0$). Consider a parallel spinor field $\psi$ on $\Delta$ (which is then a constant section of $\Sigma_0\Delta$) and let $\varphi:=P_+\psi.$ Then, the Dirac operator of $M$ associated to $\varphi$ is equal to
$$D_M\varphi=\textbf{D}_{\textbf{S}}P_+\psi=P_-\textbf{D}_{\textbf{S}}\psi=\frac{n+1}{2}HP_{-}\psi.$$
The normal bundle $Q$ of the flow is just the tangent space of $M_1$  and the connection $A^Q$ is the connection $A^M|_{M_1},$ i.e. $\theta=0.$ Since the spinor $\varphi$ is clearly basic (it is constant along the $\mathbb{S}^1$-fibers), we deduce that $D_b\varphi= \frac{n+1}{2}H\partial_t\cdot_MP_{-}\psi.$ But using the fact that $i\partial_t\cdot\psi=\psi,$ we have that
$$\partial_t\cdot_MP_{-}\psi=\partial_t\cdot\nu\cdot P_{-}\psi=\frac{i}{2}\partial_t\cdot(\psi-i\nu\cdot\psi)=P_+\psi=\varphi.$$
Therefore $\varphi$ is a solution of the basic Dirac equation with $H_0=H$ and the equality in \eqref{eq:integraline22-int} is realized.
\hfill$\square$\\\\

\noindent{\bf Acknowledgment.} Part of this work was done while Fida El Chami and Georges Habib enjoyed the hospitality of the University of Lorraine.
The second-named author would like to thank the Alexander von Humboldt Foundation for its support.

\end{document}